\documentclass[a4paper,12pt]{article} 
\usepackage{amsmath} 
\usepackage{amsopn} 
\usepackage{amssymb} 
\usepackage[mathscr]{eucal} 
\usepackage{theorem} 
\usepackage{enumerate} 

\theorembodyfont{\itshape} 
\theoremstyle{plain}
  \newtheorem{thm}{Theorem}[section] 
  \newtheorem{pro}[thm]{Proposition} 
  \newtheorem{lem}[thm]{Lemma}

\theorembodyfont{\rmfamily} 
\theoremstyle{plain}
  \newtheorem{rem}{Remark}[section] 
  \newtheorem{ex}{Example}[section] 

\renewcommand{\theequation}%
           {\thesection.\arabic{equation}}

\begin{document} 

\begin{center} 
{\LARGE The twistor lifts of surfaces in 4-spaces} 

\vspace{6mm} 

{\Large Naoya {\sc Ando}} 
\end{center} 

\vspace{3mm} 

\begin{quote} 
{\footnotesize \it Abstract} \ 
{\footnotesize 
This is a survey of the twistor lifts of surfaces in $4$-dimensional spaces. 
In most part of this survey, 
the space is Euclidean $4$-space $E^4$. 
The definitions of the Gauss maps and the twistor lifts of surfaces in $E^4$ 
are given by orthogonal complex structures of $E^4$. 
Based on these definitions, 
we can understand holomorphicity of the Gauss maps of minimal surfaces 
in $E^4$ and isotropicity of such surfaces. 
}
\end{quote} 

\vspace{3mm} 

\section{Introduction} 

\setcounter{equation}{0} 

Minimal surfaces in Euclidean $3$-space $E^3$ are closely related to 
holomorphic functions. 
The components of the immersion which gives a minimal surface in $E^3$ 
are harmonic functions. 
Considering the surface to be a Riemann surface 
by isothermal coordinates, 
we can locally construct the immersion 
by a holomorphic $1$-form and a meromorphic function 
as in the Enneper-Weierstrass representation. 
The Gauss map of a minimal surface in $E^3$ is given by 
the meromorphic function in this representation 
and it is a holomorphic map from the surface into 
the Riemann sphere $\mbox{\boldmath{$C$}} \cup \{ \infty \} 
                   =\mbox{\boldmath{$C$}}\!P^1$.  
We can refer to a well-known textbook \cite{osserman2} 
for minimal surfaces in $E^3$. 

Minimal surfaces in Euclidean $4$-space $E^4$ are also closely related to 
holomorphic functions. 
The components of the immersion giving a minimal surface in $E^4$ are 
harmonic. 
An analogue of the Enneper-Weierstrass representation 
for minimal surfaces in $E^3$ is obtained 
for minimal surfaces in $E^4$ 
(\cite{AAIK}, \cite{ando(7)}, \cite{HO}, \cite{HO2}, \cite{osserman}). 
The Gauss map of a minimal surface in $E^4$ assings to each point of 
the surface an oriented two-dimensional subspace of $E^4$ 
parallel to the tangent plane. 
Therefore it is a map from the surface into 
an oriented real Grassmann manifold 
$\Tilde{G}_{4, 2} (\mbox{\boldmath{$R$}} )$. 
In addition, 
based on identification of $\Tilde{G}_{4, 2} (\mbox{\boldmath{$R$}} )$ 
with $\mbox{\boldmath{$C$}}\!P^1 \times \mbox{\boldmath{$C$}}\!P^1$, 
it is known that the Gauss map of a minimal surface in $E^4$ is 
holomorphic (\cite{HO}). 
A holomorphic quartic differential $Q$ is defined 
on each minimal surface in $E^4$ 
(refer to \cite{calabi}, \cite{ando(3)} for such differentials). 
For minimal surfaces in $E^4$, 
we can consider the notion of isotropicity, 
and an isotropic minimal surface in $E^4$ is characterized by $Q=0$. 
It is seen that 
an isotropic minimal surface in $E^4$ is just a complex curve 
with respect to an orthogonal complex structure of $E^4$, 
which is a result of \cite{friedrich} in the special case 
where the ambient space is $E^4$.  

For a minimal surface in $E^4$, 
holomorphicity of the Gauss map and isotropicity are understood 
in terms of the twistor lifts of the surface. 
In this survey, 
the two-fold exterior power $\bigwedge^2\!E^4$ of $E^4$ is expressed 
as a $6$-dimensional subspace 
of a $16$-dimensional vector space $M(4, \mbox{\boldmath{$R$}} )$ 
(the set of $4\times 4$ matrices such that the components are real numbers). 
This subspace coincides with 
the set of alternating matrices in $M(4, \mbox{\boldmath{$R$}} )$, 
and it is decomposed into 
two $3$-dimensional subspaces $\bigwedge^2_{\pm}\!E^4$. 
Then the space of orthogonal complex structures of $E^4$ is expressed as 
the disjoint union of the unit spheres of $\bigwedge^2_{\pm}\!E^4$, 
denoted by $\Sigma_+$, $\Sigma_-$ respectively (Theorem~\ref{thm:W+-}). 
Each twistor space associated with $E^4$ is given by 
the product of $E^4$ and $\Sigma_{\varepsilon}$ for $\varepsilon =+$ or $-$. 
Therefore, in this situation, 
each twistor lift is given by a map from the surface 
into $\Sigma_{\varepsilon}$. 
In addition, 
since there exists a one-to-one correspondence 
between $\Tilde{G}_{4, 2} (\mbox{\boldmath{$R$}} )$ 
and     $\Sigma_+ \times \Sigma_-$ (Theorem~\ref{thm:GrW+-}), 
the Gauss map is given by the pair of the two twistor lifts. 
Based on the above setting in this survey, 
the twistor lifts of a surface in $E^4$ are explicitly represented 
by the components of the immersion (Proposition~\ref{pro:tl}). 
By this representation, 
we can understand holomorphicity of the Gauss map 
of a minimal surface in $E^4$ (Theorem~\ref{thm:HO}). 
In addition, 
by the holomorphic quartic differential $Q$ 
(see \eqref{Qdef} and Theorem~\ref{thm:Q} 
in Section~\ref{sect:RsurfQ}), 
we can understand isotropicity of such a surface 
(Theorem~\ref{thm:F}). 

Section~\ref{sect:notE4} is devoted to a brief overview 
of results related to Theorems~\ref{thm:HO} and \ref{thm:F} 
in $4$-dimensional spaces other than $E^4$. 

\section{Surfaces in $\mbox{\boldmath{$E^4$}}$} 

\setcounter{equation}{0} 

\subsection{The Gauss-Weingarten formulas} 

Let $O$ be an open set of the $uv$-plane $\mbox{\boldmath{$R$}}^2$ 
homeomorphic to an open disc. 
Let $F: O\longrightarrow E^4$ be 
a $C^{\infty}$-map from $O$ into Euclidean $4$-space $E^4$ 
(an $E^4$-valued $C^{\infty}$-function on $O$). 
Then $F$ is represented as  $F=(f^1 , f^2 , f^3 , f^4 )$ 
by $C^{\infty}$-functions $f^1$, $f^2$, $f^3$, $f^4$ on $O$. 
We set 
\begin{equation*} 
F_u :=(f^1_u , f^2_u , f^3_u , f^4_u ), \quad 
F_v :=(f^1_v , f^2_v , f^3_v , f^4_v ), 
\end{equation*} 
where $f^i_u$, $f^i_v$ denote the partial derivatives of $f^i$ 
with respect to the variables $u$, $v$ respectively ($i=1, 2, 3, 4$). 
Then $F_u$, $F_v$ are $E^4$-valued $C^{\infty}$-functions on $O$. 
In this survey, 
a $C^{\infty}$-map $F: O\longrightarrow E^4$ is called a \textit{surface\/} 
if at each point $p$ of $O$, 
the two vectors $F_u (p)$, $F_v (p)$ of $E^4$ are linearly independent. 

In the following, 
we suppose that $F: O\longrightarrow E^4$ is a surface. 
Let $g_{ij}$ ($i, j=1, 2$) be $C^{\infty}$-functions on $O$ defined by 
\begin{equation} 
g_{11}         :=\langle F_u , F_u \rangle , \quad 
g_{12} =g_{21} :=\langle F_u , F_v \rangle , \quad 
g_{22}         :=\langle F_v , F_v \rangle 
\label{gij} 
\end{equation} 
where $\langle \ , \ \rangle$ denotes the standard inner product 
of $E^4$, 
that is, 
for two vectors $\mbox{\boldmath{$a$}} =(a^1 , a^2 , a^3 , a^4 )$, 
                $\mbox{\boldmath{$b$}} =(b^1 , b^2 , b^3 , b^4 )$, 
we set 
\begin{equation*} 
 \langle \mbox{\boldmath{$a$}} , \mbox{\boldmath{$b$}} \rangle 
:=a^1 b^1 +a^2 b^2 +a^3 b^3 +a^4 b^4 . 
\end{equation*} 
We also set $|\mbox{\boldmath{$a$}} | 
            :=\sqrt{\langle \mbox{\boldmath{$a$}} , 
                            \mbox{\boldmath{$a$}} \rangle}$. 
Since we suppose that $F$ is a surface, 
we have $g_{11} >0$, $g_{22} >0$. 
In addition, since 
\begin{equation*} 
  g_{11} +2g_{12} t +g_{22} t^2 
=\langle F_u +tF_v , F_u +tF_v \rangle >0 
\end{equation*} 
for any $t\in \mbox{\boldmath{$R$}}$, 
we have $g_{11} g_{22} -g^2_{12} >0$. 
By $g_{ij}$ ($i, j=1, 2$), 
we can define the induced metric $g$ by $F$, 
which is a symmetric, positive-definite $2$-tensor field: 
\begin{equation*} 
g=\sum_{i, j=1, 2} g_{ij} du^i du^j \quad 
(u^1 :=u, u^2 :=v). 
\end{equation*} 

For $a, b\in \{ u, v\}$, 
we set $F_{ab} :=(f^1_{ab} , f^2_{ab} , f^3_{ab} , f^4_{ab} )$. 
We denote $F_{ab}$ by 
\begin{equation*} 
F_{ab} =F_{ab}^{\top} +F_{ab}^{\bot} , 
\end{equation*} 
where 
\begin{itemize} 
\item{$F_{ab}^{\top}$ is the tangential component of $F_{ab}$ 
and therefore represented by a linear combination of $F_u$ and $F_v$ 
at each point of $O$,}  
\item{$F_{ab}^{\bot}$ is the normal component of $F_{ab}$ 
and therefore normal to both $F_u$ and $F_v$ at each point of $O$.} 
\end{itemize} 
Let $\Gamma^k_{ij}$ ($i, j, k=1, 2$) be $C^{\infty}$ functions on $O$ 
defined by 
\begin{equation*} 
\begin{split} 
F_{uu}^{\top} & =\Gamma^1_{11} F_u +\Gamma^2_{11} F_v , \\ 
F_{uv}^{\top} & =\Gamma^1_{12} F_u +\Gamma^2_{12} F_v , \\ 
F_{vv}^{\top} & =\Gamma^1_{22} F_u +\Gamma^2_{22} F_v 
\end{split} 
\end{equation*} 
and $\Gamma^k_{21} =\Gamma^k_{12}$. 
The functions $\Gamma^k_{ij}$ ($i, j, k=1, 2$) are called 
the \textit{Christoffel symbols}. 
The following proposition is well-known: 

\begin{pro}\label{pro:Gamma} 
The Christoffel symbols $\Gamma^k_{ij}$ $(i, j, k=1, 2)$ 
are represented by $g_{11}$, $g_{12}$, $g_{22}$ 
and their partial derivatives. 
\end{pro} 

Let $\mbox{\boldmath{$n$}}_1$, $\mbox{\boldmath{$n$}}_2$ be 
$E^4$-valued $C^{\infty}$-functions on $O$ satisfying 
\begin{equation} 
\begin{split} 
& |\mbox{\boldmath{$n$}}_1 |=|\mbox{\boldmath{$n$}}_2 |=1, 
  \quad 
  \langle \mbox{\boldmath{$n$}}_1 , \mbox{\boldmath{$n$}}_2 \rangle 
 =0, \\ 
& \langle \mbox{\boldmath{$n$}}_k , F_u \rangle 
 =\langle \mbox{\boldmath{$n$}}_k , F_v \rangle 
 =0 \ \ (k=1, 2). 
\end{split} 
\label{n1n2} 
\end{equation} 
Then there exist $C^{\infty}$-functions $b^k_{ij}$ ($i, j, k=1, 2$) on $O$ 
satisfying 
\begin{equation} 
\begin{split} 
F_{uu}^{\bot} & =b^1_{11} \mbox{\boldmath{$n$}}_1 
                +b^2_{11} \mbox{\boldmath{$n$}}_2 , \\ 
F_{uv}^{\bot} & =b^1_{12} \mbox{\boldmath{$n$}}_1 
                +b^2_{12} \mbox{\boldmath{$n$}}_2 , \\ 
F_{vv}^{\bot} & =b^1_{22} \mbox{\boldmath{$n$}}_1 
                +b^2_{22} \mbox{\boldmath{$n$}}_2 
\end{split} 
\label{Fperpab} 
\end{equation} 
and $b^k_{21} =b^k_{12}$. 
We represent $\mbox{\boldmath{$n$}}_k$ 
as $\mbox{\boldmath{$n$}}_k =(n^1_k , n^2_k , n^3_k , n^4_k )$. 
In addition, we set 
\begin{equation*} 
\begin{split} 
\mbox{\boldmath{$n$}}_{ku} & :=(n^1_{ku} , n^2_{ku} , n^3_{ku} , n^4_{ku} ), 
\\ 
\mbox{\boldmath{$n$}}_{kv} & :=(n^1_{kv} , n^2_{kv} , n^3_{kv} , n^4_{kv} ), 
\end{split} 
\end{equation*} 
where $n^i_{ku}$, $n^i_{kv}$ denote the partial derivatives of $n^i_k$ 
with respect to the variables $u$, $v$ respectively ($i=1, 2, 3, 4$). 
For $a\in \{ u, v\}$, 
we denote $\mbox{\boldmath{$n$}}_{ka}$ by 
\begin{equation*} 
 \mbox{\boldmath{$n$}}_{ka} 
=\mbox{\boldmath{$n$}}_{ka}^{\top} 
+\mbox{\boldmath{$n$}}_{ka}^{\bot} , 
\end{equation*} 
where $\mbox{\boldmath{$n$}}_{ka}^{\top}$, 
      $\mbox{\boldmath{$n$}}_{ka}^{\bot}$ are 
the tangential and normal components of $\mbox{\boldmath{$n$}}_{ka}$ 
respectively. 
There exist $C^{\infty}$-functions $a^j_{kl}$ ($j, k, l=1, 2$) 
on $O$ satisfying 
\begin{equation} 
\begin{split} 
  \mbox{\boldmath{$n$}}_{ku}^{\top} 
& =-a^1_{k1} F_u -a^2_{k1} F_v , \\ 
  \mbox{\boldmath{$n$}}_{kv}^{\top} 
& =-a^1_{k2} F_u -a^2_{k2} F_v , 
\end{split} 
\label{ntopka} 
\end{equation} 
that is, 
\begin{equation} 
  (\mbox{\boldmath{$n$}}_{ku}^{\top} , \ 
   \mbox{\boldmath{$n$}}_{kv}^{\top} ) 
=-(F_u , \ F_v )A_k , \quad  
   A_k :=\left[ \begin{array}{cc} 
                 a^1_{k1} & a^1_{k2} \\ 
                 a^2_{k1} & a^2_{k2} 
                  \end{array} 
         \right] . 
\label{ntopka2} 
\end{equation} 
Therefore noticing $g_{11} g_{22} -g_{12}^2 \not= 0$, we obtain 

\begin{pro}\label{pro:A} 
The matrices $A_k$ $(k=1, 2)$ are given by 
\begin{equation*} 
A_k = 
\left[ 
\begin{array}{cc} 
 g_{11} & g_{12} \\ 
 g_{21} & g_{22}  
  \end{array} 
\right]^{-1}  
\left[ 
\begin{array}{cc} 
 b^k_{11} & b^k_{12} \\ 
 b^k_{21} & b^k_{22} 
  \end{array} 
\right] . 
\end{equation*} 
\end{pro} 

\begin{rem} 
We see that $A_k$ is the representation matrix of the shape operator $S_k$ 
for the normal vector field $\mbox{\boldmath{$n$}}_k$, 
which is a $(1, 1)$-tensor field symmetric with respect to 
the induced metric $g$: 
\begin{equation*} 
(S_k (F_u ), S_k (F_v )) 
=(F_u , F_v )A_k . 
\end{equation*} 
\end{rem} 

Since $|\mbox{\boldmath{$n$}}_1 |=1$, 
we have $\langle \mbox{\boldmath{$n$}}_{1a} , 
                 \mbox{\boldmath{$n$}}_1 \rangle =0$ for $a=u, v$. 
Therefore there exist functions $\gamma_1$, $\gamma_2$ satisfying 
\begin{equation*} 
\mbox{\boldmath{$n$}}^{\bot}_{1u} =\gamma_1 \mbox{\boldmath{$n$}}_2 , \quad 
\mbox{\boldmath{$n$}}^{\bot}_{1v} =\gamma_2 \mbox{\boldmath{$n$}}_2 . 
\end{equation*} 
Then noticing $\langle \mbox{\boldmath{$n$}}_1 , 
                       \mbox{\boldmath{$n$}}_2 \rangle =0$ 
and                  $|\mbox{\boldmath{$n$}}_2 |=1$, 
we have 
\begin{equation*} 
\mbox{\boldmath{$n$}}^{\bot}_{2u} =-\gamma_1 \mbox{\boldmath{$n$}}_1 , \quad 
\mbox{\boldmath{$n$}}^{\bot}_{2v} =-\gamma_2 \mbox{\boldmath{$n$}}_1 . 
\end{equation*} 
Hence we obtain 

\begin{thm}\label{thm:GW} 
A surface $F$ satisfies the Gauss formulas 
\begin{equation*} 
\begin{split} 
  F_{uu} 
& =\Gamma^1_{11} F_u +\Gamma^2_{11} F_v 
  + b^1_{11} \mbox{\boldmath{$n$}}_1 
  +b^2_{11} \mbox{\boldmath{$n$}}_2 , \\ 
  F_{uv} 
& =\Gamma^1_{12} F_u +\Gamma^2_{12} F_v 
  + b^1_{12} \mbox{\boldmath{$n$}}_1 
  + b^2_{12} \mbox{\boldmath{$n$}}_2 , \\ 
  F_{vv} 
& =\Gamma^1_{22} F_u +\Gamma^2_{22} F_v 
  + b^1_{22} \mbox{\boldmath{$n$}}_1 
  + b^2_{22} \mbox{\boldmath{$n$}}_2 
\end{split} 
\end{equation*} 
and the Weingarten formulas 
\begin{equation*} 
\begin{split} 
  \mbox{\boldmath{$n$}}_{1u} 
& =-a^1_{11} F_u -a^2_{11} F_v 
   +\gamma_1 \mbox{\boldmath{$n$}}_2 , \\ 
  \mbox{\boldmath{$n$}}_{1v} 
& =-a^1_{12} F_u -a^2_{12} F_v 
   +\gamma_2 \mbox{\boldmath{$n$}}_2 , \\ 
  \mbox{\boldmath{$n$}}_{2u} 
& =-a^1_{21} F_u -a^2_{21} F_v 
   -\gamma_1 \mbox{\boldmath{$n$}}_1 , \\ 
  \mbox{\boldmath{$n$}}_{2v} 
& =-a^1_{22} F_u -a^2_{22} F_v 
   -\gamma_2 \mbox{\boldmath{$n$}}_1 , 
\end{split} 
\end{equation*} 
that is, 
\begin{equation*} 
\begin{split} 
   (F_u , \ F_v , \ \mbox{\boldmath{$n$}}_1 , \ 
                    \mbox{\boldmath{$n$}}_2 )_u 
& =(F_u , \ F_v , \ \mbox{\boldmath{$n$}}_1 , \ 
                    \mbox{\boldmath{$n$}}_2 )R_1 , \\ 
   (F_u , \ F_v , \ \mbox{\boldmath{$n$}}_1 , \ 
                    \mbox{\boldmath{$n$}}_2 )_v 
& =(F_u , \ F_v , \ \mbox{\boldmath{$n$}}_1 , \ 
                    \mbox{\boldmath{$n$}}_2 )R_2 , 
\end{split} 
\end{equation*} 
where 
\begin{equation*} 
R_j =\left[ 
     \begin{array}{cccc} 
      \Gamma^1_{1j} & \Gamma^1_{2j} & - a^1_{1j} & - a^1_{2j} \\ 
      \Gamma^2_{1j} & \Gamma^2_{2j} & - a^2_{1j} & - a^2_{2j} \\ 
           b^1_{1j} &      b^1_{2j} &   0        & -\gamma_j  \\ 
           b^2_{1j} &      b^2_{2j} &  \gamma_j  &   0 
    \end{array} 
  \right] \ (j=1, 2). 
\end{equation*} 
\end{thm}

\subsection{The mean curvature vectors of surfaces 
in $\mbox{\boldmath{$E^4$}}$}\label{subsect:H} 

Let $\Tilde{\mbox{\boldmath{$n$}}}_1$, 
    $\Tilde{\mbox{\boldmath{$n$}}}_2$ be 
$E^4$-valued $C^{\infty}$-functions on $O$ 
satisfying 
\begin{equation*} 
 |\Tilde{\mbox{\boldmath{$n$}}}_1 | 
=|\Tilde{\mbox{\boldmath{$n$}}}_2 | 
=1, 
 \ 
 \langle \Tilde{\mbox{\boldmath{$n$}}}_1 , 
         \Tilde{\mbox{\boldmath{$n$}}}_2 \rangle 
=0, 
 \ 
 \langle \Tilde{\mbox{\boldmath{$n$}}}_k , F_u \rangle 
=\langle \Tilde{\mbox{\boldmath{$n$}}}_k , F_v \rangle 
=0 \quad (k=1, 2).  
\end{equation*} 
Then there exists a $C^{\infty}$-function $\theta$ on $O$ 
satisfying either 
\begin{equation} 
 (\Tilde{\mbox{\boldmath{$n$}}}_1 , \ \Tilde{\mbox{\boldmath{$n$}}}_2 ) 
=(       \mbox{\boldmath{$n$}}_1  , \        \mbox{\boldmath{$n$}}_2  )
  \left[ \begin{array}{cc} 
          \cos \theta & -\sin \theta \\ 
          \sin \theta &  \cos \theta 
           \end{array} 
  \right] 
\label{+} 
\end{equation} 
or 
\begin{equation} 
 (\Tilde{\mbox{\boldmath{$n$}}}_1 , \ \Tilde{\mbox{\boldmath{$n$}}}_2 ) 
=(       \mbox{\boldmath{$n$}}_1  , \        \mbox{\boldmath{$n$}}_2  )
  \left[ \begin{array}{cc} 
          \cos \theta &  \sin \theta \\ 
          \sin \theta & -\cos \theta 
           \end{array} 
  \right] . 
\label{-} 
\end{equation} 

Suppose \eqref{+}. Then for $a\in \{ u, v\}$, we obtain 
\begin{equation} 
\begin{split} 
    \Tilde{\mbox{\boldmath{$n$}}}^{\top}_{1a} 
& =((\cos \theta )\mbox{\boldmath{$n$}}_1 
   +(\sin \theta )\mbox{\boldmath{$n$}}_2 )^{\top}_a \\ 
& =(-\theta_a (\sin \theta )\mbox{\boldmath{$n$}}_1 
    +         (\cos \theta )\mbox{\boldmath{$n$}}_{1a} 
    +\theta_a (\cos \theta )\mbox{\boldmath{$n$}}_2 
    +         (\sin \theta )\mbox{\boldmath{$n$}}_{2a} )^{\top} \\ 
& = (\cos \theta )\mbox{\boldmath{$n$}}^{\top}_{1a} 
   +(\sin \theta )\mbox{\boldmath{$n$}}^{\top}_{2a} . 
\end{split} 
\label{Tildentop1a0} 
\end{equation} 
Similarly, we obtain 
\begin{equation} 
    \Tilde{\mbox{\boldmath{$n$}}}^{\top}_{2a} 
 =-(\sin \theta )\mbox{\boldmath{$n$}}^{\top}_{1a} 
  +(\cos \theta )\mbox{\boldmath{$n$}}^{\top}_{2a} . 
\label{Tildentop2a} 
\end{equation} 
Let $\Tilde{A}_k$ ($k=1, 2$) be matrices defined by 
\begin{equation} 
  (\Tilde{\mbox{\boldmath{$n$}}}_{ku}^{\top} , \ 
   \Tilde{\mbox{\boldmath{$n$}}}_{kv}^{\top} ) 
=-(F_u , \ F_v )\Tilde{A}_k . 
\label{tildeAk} 
\end{equation} 

\begin{lem}\label{lem:TildeA} 
Suppose \eqref{+}. 
Then the matrices $\Tilde{A}_k$ $(k=1, 2)$ are represented as 
\begin{equation*} 
\Tilde{A}_1 = (\cos \theta )A_1 +(\sin \theta )A_2 , \quad 
\Tilde{A}_2 =-(\sin \theta )A_1 +(\cos \theta )A_2 . 
\end{equation*} 
\end{lem} 

\vspace{3mm} 

\par\noindent 
\textit{Proof} \ 
By \eqref{Tildentop1a0}, we already know 
\begin{equation*} 
  \Tilde{\mbox{\boldmath{$n$}}}^{\top}_{1a} 
=(\cos \theta )\mbox{\boldmath{$n$}}^{\top}_{1a} 
+(\sin \theta )\mbox{\boldmath{$n$}}^{\top}_{2a} 
\end{equation*} 
for $a=u, v$. 
Combining this with \eqref{ntopka}, we obtain 
\begin{equation*} 
\begin{split} 
& (\Tilde{\mbox{\boldmath{$n$}}}_{1u}^{\top} , \ 
   \Tilde{\mbox{\boldmath{$n$}}}_{1v}^{\top} ) \\ 
& =-(F_u , \ F_v )
   \left[ \begin{array}{cc} 
           (\cos \theta )a^1_{11} +(\sin \theta )a^1_{21} & 
           (\cos \theta )a^1_{12} +(\sin \theta )a^1_{22} \\ 
           (\cos \theta )a^2_{11} +(\sin \theta )a^2_{21} & 
           (\cos \theta )a^2_{12} +(\sin \theta )a^2_{22} 
            \end{array} 
   \right] \\ 
& =-(F_u , \ F_v )((\cos \theta )A_1 +(\sin \theta )A_2 ). 
\end{split} 
\end{equation*} 
This means 
\begin{equation} 
\Tilde{A}_1 = (\cos \theta )A_1 +(\sin \theta )A_2 . 
\label{tildeA1} 
\end{equation} 
Similarly, combining \eqref{ntopka} with \eqref{Tildentop2a}, we obtain 
\begin{equation*} 
  (\Tilde{\mbox{\boldmath{$n$}}}_{2u}^{\top} , \ 
   \Tilde{\mbox{\boldmath{$n$}}}_{2v}^{\top} ) 
=-(F_u , \ F_v )(-(\sin \theta )A_1 +(\cos \theta )A_2 ). 
\end{equation*} 
This means 
\begin{equation*} 
\Tilde{A}_2 =-(\sin \theta )A_1 +(\cos \theta )A_2 . 
\end{equation*}  
Hence we have proved Lemma~\ref{lem:TildeA}. 
\hfill 
$\square$ 

\vspace{3mm} 

In the above discussion, we supposed \eqref{+}. 
If we suppose \eqref{-} instead of \eqref{+}, 
then by a similar discussion, we obtain 

\begin{lem}\label{lem:TildeA2} 
Suppose \eqref{-}. 
Then the matrices $\Tilde{A}_k$ $(k=1, 2)$ are represented as 
\begin{equation*} 
\Tilde{A}_1 = (\cos \theta )A_1 +(\sin \theta )A_2 , \quad 
\Tilde{A}_2 = (\sin \theta )A_1 -(\cos \theta )A_2 . 
\end{equation*} 
\end{lem} 

\begin{rem}\label{rem:rmso}  
From \eqref{tildeAk} and \eqref{tildeA1}, 
we see that $(\cos \theta )A_1 +(\sin \theta )A_2$ is 
the representation matrix of the shape operator 
for the normal vector field $\mbox{\boldmath{$n$}} (\theta ) 
             :=(\cos \theta )\mbox{\boldmath{$n$}}_1 
              +(\sin \theta )\mbox{\boldmath{$n$}}_2$ with length one. 
\end{rem} 

We will prove 

\begin{pro}\label{pro:H} 
An $E^4$-valued $C^{\infty}$-function on $O$ 
defined by 
\begin{equation*} 
  \mbox{\boldmath{$H$}} 
:=\dfrac{1}{2} (({\rm tr}\,A_1 )\mbox{\boldmath{$n$}}_1 
               +({\rm tr}\,A_2 )\mbox{\boldmath{$n$}}_2 )
\end{equation*} 
does not depend on the choice of 
a pair $(\mbox{\boldmath{$n$}}_1 , \mbox{\boldmath{$n$}}_2 )$ 
satisfying \eqref{n1n2}. 
\end{pro} 

\vspace{3mm} 

\par\noindent 
\textit{Proof} \ 
Let $\Tilde{\mbox{\boldmath{$n$}}}_1$, 
    $\Tilde{\mbox{\boldmath{$n$}}}_2$ and $\theta$ be as above. 
    Then we have Lemma~\ref{lem:TildeA} or \ref{lem:TildeA2}. 
If we have Lemma~\ref{lem:TildeA}, then we obtain 
\begin{equation*} 
\begin{split} 
& ({\rm tr}\,\Tilde{A}_1 )\Tilde{\mbox{\boldmath{$n$}}}_1 
 +({\rm tr}\,\Tilde{A}_2 )\Tilde{\mbox{\boldmath{$n$}}}_2 \\ 
& =( (\cos \theta )\,{\rm tr}\,A_1 +(\sin \theta )\,{\rm tr}\,A_2 ) 
   ( (\cos \theta )\mbox{\boldmath{$n$}}_1 
    +(\sin \theta )\mbox{\boldmath{$n$}}_2 ) \\ 
&    \quad \quad 
  +(-(\sin \theta )\,{\rm tr}\,A_1 +(\cos \theta )\,{\rm tr}\,A_2 ) 
   (-(\sin \theta )\mbox{\boldmath{$n$}}_1 
    +(\cos \theta )\mbox{\boldmath{$n$}}_2 ) \\ 
& =(\cos^2 \theta +\sin^2 \theta )({\rm tr}\,A_1 )\mbox{\boldmath{$n$}}_1 
  +(\cos^2 \theta +\sin^2 \theta )({\rm tr}\,A_2 )\mbox{\boldmath{$n$}}_2 \\ 
& =({\rm tr}\,A_1 )\mbox{\boldmath{$n$}}_1 
  +({\rm tr}\,A_2 )\mbox{\boldmath{$n$}}_2 . 
\end{split} 
\end{equation*} 
If we have Lemma~\ref{lem:TildeA2}, 
then we similarly obtain the same result. 
Hence we obtain Proposition~\ref{pro:H}. 
\hfill 
$\square$ 

\vspace{3mm} 

The $E^4$-valued $C^{\infty}$-function 
$\mbox{\boldmath{$H$}}$ is called 
the \textit{mean curvature vector\/} (\textit{field\/}) 
of a surface $F$. 
A surface $F$ is said to be \textit{minimal\/} 
if the mean curvature vector $\mbox{\boldmath{$H$}}$ of $F$ is 
identically zero on $O$. 

In the following, 
we suppose that $(u, v)$ are \textit{isothermal\/} coordinates 
of the induced metric $g$, 
that is, 
$(u, v)$ satisfy $g_{11} =g_{22}$ and $g_{12} =0$. 
Let $\alpha$ be a $C^{\infty}$-function on $O$ 
satisfying $e^{2\alpha} =g_{11}$. 

\begin{pro}\label{pro:ic}
The following hold\/$:$ 
\begin{equation*} 
\begin{split} 
& \Gamma^1_{11} =\Gamma^2_{12} =-\Gamma^1_{22} =\alpha_u , \quad 
 -\Gamma^2_{11} =\Gamma^1_{12} = \Gamma^2_{22} =\alpha_v , \\ 
& a^j_{kl} =\dfrac{1}{e^{2\alpha}} b^k_{jl} , \quad 
  \mbox{\boldmath{$H$}} 
 =\dfrac{1}{2e^{2\alpha}} 
  \sum^2_{k=1} (b^k_{11} +b^k_{22} )\mbox{\boldmath{$n$}}_k . 
\end{split} 
\end{equation*} 
\end{pro} 

Suppose that $F$ is minimal. 
Then from Theorem~\ref{thm:GW} and Proposition~\ref{pro:ic}, 
we obtain $F_{uu} +F_{vv} =(0, 0, 0, 0)$. 
Since we represent $F$ as $F=(f^1 , f^2 , f^3 , f^4 )$, 
each $f^i$ is a harmonic function, that is, 
$f=f^i$ is a solution of the Laplace 
equation $f_{uu} +f_{vv} =0$. 
If we set 
\begin{equation} 
\beta^k :=\dfrac{1}{2} (b^k_{11} -\sqrt{-1} b^k_{12} ) \ \ 
(k=1, 2), \quad 
\gamma  :=\dfrac{1}{2} (\gamma_1 +\sqrt{-1} \gamma_2 ) 
\label{betagamma} 
\end{equation} 
and 
\begin{equation} 
  \dfrac{\partial}{\partial w} 
:=\dfrac{1}{2} \left( \dfrac{\partial}{\partial u} 
                     -\sqrt{-1} 
                      \dfrac{\partial}{\partial v} \right) , \quad 
  \dfrac{\partial}{\partial \overline{w}} 
:=\dfrac{1}{2} \left( \dfrac{\partial}{\partial u} 
                     +\sqrt{-1} 
                      \dfrac{\partial}{\partial v} \right) , 
\label{dw} 
\end{equation} 
then computing the both sides of 
\begin{equation*} 
F_{uuv} =F_{uvu} , \quad 
F_{uvv} =F_{vvu} , \quad 
\mbox{\boldmath{$n$}}_{kuv} =\mbox{\boldmath{$n$}}_{kvu} \ (k=1, 2) 
\end{equation*} 
or 
\begin{equation*} 
(R_1 )_v -(R_2 )_u =R_1 R_2 -R_2 R_1 
\end{equation*} 
by Theorem~\ref{thm:GW}, 
we obtain 

\begin{thm}\label{thm:GCR} 
For a minimal surface $F$, 
the equation of Gauss is represented as 
\begin{equation*} 
 \alpha_{uu} +\alpha_{vv} 
=\dfrac{4}{e^{2\alpha}} (|\beta^1 |^2 +|\beta^2 |^2 ), 
\end{equation*} 
the equations of Codazzi are represented as 
\begin{equation*} 
\dfrac{\partial \beta^1}{\partial \overline{w}} = \beta^2 \gamma , \quad 
\dfrac{\partial \beta^2}{\partial \overline{w}} =-\beta^1 \gamma 
\end{equation*} 
and the equation of Ricci is represented as 
\begin{equation*} 
{\rm Im}\,\left( \dfrac{\partial \gamma}{\partial w} 
                +\dfrac{2}{e^{2\alpha}} \beta^1 \overline{\beta}^2 
           \right) =0. 
\end{equation*} 
\end{thm} 

By the equations of Codazzi in Theorem~\ref{thm:GCR}, 
we obtain 

\begin{pro}\label{pro:b1b2} 
For a minimal surface $F$, 
a complex-valued function $(\beta^1 )^2 +(\beta^2 )^2$ is holomorphic, 
that is, $(\beta^1 )^2 +(\beta^2 )^2$ satisfies 
$\partial ((\beta^1 )^2 +(\beta^2 )^2 )/\partial \overline{w} =0$. 
\end{pro} 

\section{Riemann surfaces and 
complex quartic differentials}\label{sect:RsurfQ} 

\setcounter{equation}{0} 

Let $M$ be a Hausdorff space. 
Then $M$ is called a \textit{two-dimensional topological manifold\/} 
if for each point $p$ of $M$, 
there exists a neighborhood $U$ of $p$ 
which is homeomorphic to an open set of $\mbox{\boldmath{$R$}}^2$. 
Let $M$ be a two-dimensional topological manifold. 
let $\{ U_{\lambda} \}_{\lambda \in \Lambda}$ be an open cover of $M$ 
such that for each $\lambda \in \Lambda$, 
$U_{\lambda}$ is homeomorphic to an open set $O_{\lambda}$ 
of $\mbox{\boldmath{$C$}} \cong \mbox{\boldmath{$R$}}^2$ 
by a homeomorphism $w_{\lambda} :U_{\lambda} \longrightarrow O_{\lambda}$. 
Then $\{ (U_{\lambda} , w_{\lambda} )\}_{\lambda \in \Lambda}$ is called 
a \textit{system of holomorphic coordinate neighborhoods\/} of $M$ 
if for $\lambda , \mu \in \Lambda$ 
with $U_{\lambda} \cap U_{\mu} \not= \emptyset$, 
a homeomorphism 
\begin{equation*} 
w_{\mu} \circ w_{\lambda}^{-1} : 
w_{\lambda} (U_{\lambda} \cap U_{\mu} )\longrightarrow 
w_{\mu}  (U_{\lambda} \cap U_{\mu} )
\end{equation*} 
from an open set $w_{\lambda} (U_{\lambda} \cap U_{\mu} )$ 
of $\mbox{\boldmath{$C$}}$ 
onto an open set $w_{\mu}  (U_{\lambda} \cap U_{\mu} )$ 
of $\mbox{\boldmath{$C$}}$ is holomorphic. 
Let $\{ (U_{\lambda} , w_{\lambda} )\}_{\lambda \in \Lambda}$ be 
a system of holomorphic coordinate neighborhoods of $M$. 
Then $M$ equipped with $\{ (U_{\lambda} , w_{\lambda} ) 
                        \}_{\lambda \in \Lambda}$ is 
called a \textit{Riemann surface}. 
For an open set $U$ and a homeomorphism $w$ from $U$ 
onto an open set $O$ of $\mbox{\boldmath{$C$}}$, 
if $w_{\lambda} \circ w^{-1}$ is a homeomorphism 
from $w(U_{\lambda} \cap U)$ onto $w_{\lambda} (U_{\lambda} \cap U)$ 
whenever $U_{\lambda} \cap U\not= \emptyset$, 
then the pair $(U, w)$ is called 
a \textit{holomorphic coordinate neighborhood\/} of $M$. 
Therefore, each $(U_{\lambda} , w_{\lambda} )$ is 
a holomorphic coordinate neighborhood of $M$. 
For a holomorphic coordinate neighborhood $(U, w)$ of $M$, 
we call $w$ a \textit{local complex coordinate\/} of $M$. 
We represent each local complex coordinate $w_{\lambda}$ 
as $w_{\lambda} =u_{\lambda} +\sqrt{-1} v_{\lambda}$. 
Then $\{ (U_{\lambda} , (u_{\lambda} , v_{\lambda} )) 
      \}_{\lambda \in \Lambda}$ is 
a system of real-analytic coordinate neighborhoods of $M$. 
Since $w_{\mu} \circ w_{\lambda}^{-1}$ is holomorphic, 
we have the Cauchy-Riemann equations: 
\begin{equation*} 
  \dfrac{\partial u_{\mu}}{\partial u_{\lambda}} 
= \dfrac{\partial v_{\mu}}{\partial v_{\lambda}} , \quad 
  \dfrac{\partial v_{\mu}}{\partial u_{\lambda}} 
=-\dfrac{\partial u_{\mu}}{\partial v_{\lambda}} . 
\end{equation*} 
Therefore the Jacobian $\partial (u_{\mu}  , v_{\mu}  ) 
                       /\partial (u_{\lambda} , v_{\lambda} )$ is positive, 
which means that $M$ is orientable. 

\begin{ex}\label{ex:Rsphere} 
In this survey, 
it is important that the unit $2$-sphere $S^2$ is naturally considered 
to be a special Riemann surface, i.e., 
the Riemann sphere $\mbox{\boldmath{$C$}} \cup \{ \infty \} 
                   =\mbox{\boldmath{$C$}}\!P^1$, 
via stereographic projections. 
For the unit $2$-sphere 
\begin{equation*} 
S^2 =\{ (x, y, z)\in \mbox{\boldmath{$R$}}^3 \ | \ 
         x^2 +y^2 +z^2 =1\} , 
\end{equation*} 
we set $p_+ :=(0, 0, 1)$, $p_- :=(0, 0, -1)$. 
Let $\pi_+$, $\pi_-$ be the stereographic projections 
from $p_+$, $p_-$ respectively. 
Then $\pi_{\pm}$ are the bijective maps 
from $U_{\pm} :=S^2 \setminus \{ p_{\pm} \}$ 
onto $\mbox{\boldmath{$R$}}^2$ defined by 
\begin{equation} 
  \pi_{\pm} (x, y, z)
:=\left( \dfrac{x}{1\mp z} , \ 
         \dfrac{y}{1\mp z} 
  \right) 
\label{pipm} 
\end{equation} 
for $(x, y, z)\in U_{\pm}$ respectively. 
For $\varepsilon \in \{ +, -\}$ and $(x, y, z)\in U_{\varepsilon}$, 
the intersection of the straight line 
through $(x, y, z)$ and $p_{\varepsilon}$ with the $xy$-plane 
in the $xyz$-space $\mbox{\boldmath{$R$}}^3$ is given 
by $\pi_{\pm} (x, y, z)$. 
We set 
\begin{equation*} 
w_{\pm} :=\dfrac{1}{1\mp z} (x\pm \sqrt{-1} y) \quad 
((x, y, z)\in U_{\pm} ). 
\end{equation*} 
Then for $(x, y, z)\in U_+ \cap U_-$, 
$w_+$ and $w_-$ are nonzero and satisfy $w_+ w_- =1$. 
Therefore $\{ (U_+ , \pi_+ ), (U_- , \pi_- )\}$ is 
a system of holomorphic coordinate neighborhoods of $S^2$. 
Hence $S^2$ is considered to be a Riemann surface, 
which is called the \textit{Riemann sphere}. 
The Riemann sphere can be identified with 
the one-dimensional complex projective space 
$\mbox{\boldmath{$C$}}\!P^1 
=\mbox{\boldmath{$C$}} \cup \{ \infty \}$. 
\end{ex} 

Let $M$ be a Riemann surface 
with a system $\{ (U_{\lambda} , w_{\lambda} )\}_{\lambda \in \Lambda}$ 
of holomorphic coordinate neighborhoods. 
Let $T^* M$ be the cotangent bundle of $M$ 
and let $\otimes^k T^* M$ denote 
the tensor product of $k$-copies of $T^* M$. 
Therefore the fiber of $\otimes^k T^* M$ on each point $p$ of $M$ 
consists of the covariant $k$-tensors on $T_p M$. 
A \textit{complex $k$-differential\/} on $M$ is a section of 
the complexification $\otimes^k T^* M\otimes \mbox{\boldmath{$C$}}$ 
of $\otimes^k T^* M$ which is 
locally represented as $f_{\lambda} dw_{\lambda}^k$ 
by a complex-valued function $f_{\lambda}$ on $U_{\lambda}$ 
and the tensor product $dw_{\lambda}^k$ of $k$-copies of 
the differential $dw_{\lambda}$ of $w_{\lambda}$. 
A complex $k$-differential is said to be \textit{holomorphic\/} 
if the function $f_{\lambda}$ 
on each holomorphic coordinate neighborhood $(U_{\lambda} , w_{\lambda} )$ 
is holomorphic. 
Notice that whether $f_{\lambda}$ is holomorphic or not is determined by 
the given complex $k$-differential: 
for $\lambda$, $\mu \in \Lambda$ 
with $U_{\lambda} \cap U_{\mu} \not= \emptyset$, 
if $f_{\lambda}$ is holomorphic on $U_{\lambda} \cap U_{\mu}$, 
then $f_{\mu}$ is also holomorphic, 
because $f_{\mu}$ is represented as 
\begin{equation*} 
 f_{\mu} 
=f_{\lambda} \left( \dfrac{dw_{\lambda}}{dw_{\mu}} \right)^k 
\end{equation*} 
and $w_{\lambda} \circ w_{\mu}^{-1}$ is holomorphic. 

A complex $2$-differential is also called 
a \textit{complex quadratic differential}. 
For example, 
the Hopf differential on a surface in $E^3$ is 
a complex quadratic differential. 
If the surface has constant mean curvature, 
then by the equations of Codazzi-Mainardi, 
the Hopf differential is holomorphic. 

A complex $4$-differential is also called 
a \textit{complex quartic differential}. 
Complex quartic differentials on surfaces in $E^4$ are defined 
as analogues of the Hopf differentials on surfaces in $E^3$. 
Let $M$, $\{ (U_{\lambda} , w_{\lambda} )\}_{\lambda \in \Lambda}$ be 
as above. 
Let $F$ be an immersion of $M$ into $E^4$. 
Then $F$ is said to be \textit{conformal\/} 
if $(u_{\lambda} , v_{\lambda} )$ are isothermal coordinates 
of the induced metric $g$ by $F$ on $U_{\lambda}$ 
for each $\lambda \in \Lambda$. 
Notice that if $(u_{\lambda} , v_{\lambda} )$ are isothermal coordinates, 
then for $\mu \in \Lambda$ with $U_{\lambda} \cap U_{\mu} \not= \emptyset$, 
$(u_{\mu} , v_{\mu} )$ are isothermal coordinates 
on $U_{\lambda} \cap U_{\mu}$. 
Let $F:M\longrightarrow E^4$ be a conformal immersion 
and let $dF$ denote the differential of $F$. 
Then for each $p\in M$, 
$dF$ gives a linear map $dF_p$ from $T_p M$ into $T_{F(p)} E^4$. 
In addition, $dF$ naturally gives a complex linear map 
from $T_p M\otimes \mbox{\boldmath{$C$}}$ 
into $T_{F(p)} E^4 \otimes \mbox{\boldmath{$C$}}$, 
which is also denoted by $dF_p$. 
It is obvious that $T_{F(p)} E^4 \otimes \mbox{\boldmath{$C$}}$ can be 
identified with $E^4 \otimes \mbox{\boldmath{$C$}}$ 
by parallelism in $E^4$. 
Then on each $U_{\lambda}$, 
an $E^4 \otimes \mbox{\boldmath{$C$}}$-valued function $F_{w_{\lambda}}$ 
is defined by 
\begin{equation*} 
   F_{w_{\lambda}} (p) 
:=dF_p \left( \dfrac{\partial}{\partial w_{\lambda}} \right) 
\end{equation*} 
for each $p\in U_{\lambda}$, 
where $\partial /\partial w_{\lambda}$ is defined as in \eqref{dw}. 
Since $F$ is conformal, 
we have $\langle F_{w_{\lambda}} , F_{w_{\lambda}} \rangle =0$ 
on $U_{\lambda}$, 
by considering the inner product $\langle \ , \ \rangle$ of $E^4$ 
to be a complex bilinear form on $E^4 \otimes \mbox{\boldmath{$C$}}$. 
Let $F_{w_{\lambda} w_{\lambda}}$ be 
an $E^4 \otimes \mbox{\boldmath{$C$}}$-valued function on $U_{\lambda}$ 
defined by 
\begin{equation*} 
   F_{w_{\lambda} w_{\lambda}} (p) 
:=\dfrac{\partial F_{w_{\lambda}}}{\partial w_{\lambda}} (p) 
  \quad 
(p\in U_{\lambda} ). 
\end{equation*} 
Let ${\rm Re}\,F_{w_{\lambda} w_{\lambda}}$, 
    ${\rm Im}\,F_{w_{\lambda} w_{\lambda}}$ 
be the real part and the imaginary part of $F_{w_{\lambda} w_{\lambda}}$ 
respectively. 
Then they are represented as 
\begin{equation*} 
\begin{split} 
    {\rm Re}\,F_{w_{\lambda} w_{\lambda}} 
& = {\rm Re}\,F_{w_{\lambda} w_{\lambda}}^{\top}  
   +{\rm Re}\,F_{w_{\lambda} w_{\lambda}}^{\bot} , \\  
    {\rm Im}\,F_{w_{\lambda} w_{\lambda}} 
& = {\rm Im}\,F_{w_{\lambda} w_{\lambda}}^{\top}  
   +{\rm Im}\,F_{w_{\lambda} w_{\lambda}}^{\bot} , 
\end{split} 
\end{equation*} 
where ${\rm Re}\,F_{w_{\lambda} w_{\lambda}}^{\top}$, 
      ${\rm Im}\,F_{w_{\lambda} w_{\lambda}}^{\top}$ 
are the tangential components 
and   ${\rm Re}\,F_{w_{\lambda} w_{\lambda}}^{\bot}$, 
      ${\rm Im}\,F_{w_{\lambda} w_{\lambda}}^{\bot}$ 
are the normal components 
of ${\rm Re}\,F_{w_{\lambda} w_{\lambda}}$, 
   ${\rm Im}\,F_{w_{\lambda} w_{\lambda}}$ respectively 
with respect to $F$. 
We set 
\begin{equation*} 
            F_{w_{\lambda} w_{\lambda}}^{\bot} 
:={\rm Re}\,F_{w_{\lambda} w_{\lambda}}^{\bot} 
  +\sqrt{-1} 
  {\rm Im}\,F_{w_{\lambda} w_{\lambda}}^{\bot} . 
\end{equation*} 
Then for a normal vector field $\mbox{\boldmath{$n$}}$ of $F$ 
on $U_{\lambda}$, 
we have 
\begin{equation} 
  \langle F_{w_{\lambda} w_{\lambda}}^{\bot} , \mbox{\boldmath{$n$}} \rangle 
=-\langle F_{w_{\lambda}} , \mbox{\boldmath{$n$}}_{w_{\lambda}} \rangle , 
\label{Fww} 
\end{equation} 
where $\mbox{\boldmath{$n$}}_{w_{\lambda}}$ is 
an $E^4 \otimes \mbox{\boldmath{$C$}}$-valued function on $U_{\lambda}$ 
defined by 
\begin{equation*} 
\mbox{\boldmath{$n$}}_{w_{\lambda}} (p) 
:=\dfrac{\partial \mbox{\boldmath{$n$}}}{\partial w_{\lambda}} (p) 
  \quad 
(p\in U_{\lambda} ). 
\end{equation*} 
By \eqref{Fww}, 
for $\mu \in \Lambda$ with $U_{\lambda} \cap U_{\mu} \not= \emptyset$, 
we have 
\begin{equation} 
  F_{w_{\mu} w_{\mu}}^{\bot} 
=\left( \dfrac{dw_{\lambda}}{dw_{\mu}} \right)^2 
  F_{w_{\lambda} w_{\lambda}}^{\bot} 
\label{Fab} 
\end{equation} 
on $U_{\lambda} \cap U_{\mu}$. 
We set 
\begin{equation*} 
   \Phi_{\lambda} 
:= \langle F_{w_{\lambda} w_{\lambda}}^{\bot} , 
           F_{w_{\lambda} w_{\lambda}}^{\bot} 
   \rangle . 
\end{equation*} 
Then by \eqref{Fab}, we have 
\begin{equation*} 
 \Phi_{\mu} 
=\left( \dfrac{dw_{\lambda}}{dw_{\mu}} \right)^4 
 \Phi_{\lambda} . 
\end{equation*} 
Therefore we can define a complex quartic differential $Q$ on $M$ 
by setting 
\begin{equation} 
Q=\Phi dw^4 , \quad 
\Phi =\langle F_{ww}^{\bot} , F_{ww}^{\bot} \rangle 
\label{Qdef} 
\end{equation} 
on each holomorphic coordinate neighborhood $(U, w)$ of $M$. 

\begin{rem} 
Let $M$ be a Riemann surface 
with a system $\{ (U_{\lambda} , w_{\lambda} )\}_{\lambda \in \Lambda}$ 
of holomorphic coordinate neighborhoods. 
Let $F:M\longrightarrow E^4$ be a conformal immersion. 
According to discussions in Subsection~\ref{subsect:H}, 
the mean curvature vector of $F$ is defined on 
each coordinate neighborhood $\{ (U_{\lambda} , (u_{\lambda} , v_{\lambda} )) 
                              \}_{\lambda \in \Lambda}$ 
and therefore denoted by $\mbox{\boldmath{$H$}}_{\lambda}$. 
In addition, 
noticing \eqref{ntopka2}, we see that 
if $\lambda$, $\mu \in \Lambda$ satisfy 
$U_{\lambda} \cap U_{\mu} \not= \emptyset$, 
then $\mbox{\boldmath{$H$}}_{\lambda}$ coincides with 
     $\mbox{\boldmath{$H$}}_{\mu}$ 
on $U_{\lambda} \cap U_{\mu}$. 
Hence we obtain the mean curvature vector field $\mbox{\boldmath{$H$}}$ 
of $F$ defined on $M$ 
by setting $\mbox{\boldmath{$H$}} =\mbox{\boldmath{$H$}}_{\lambda}$ 
on each $U_{\lambda}$. 
In particular, 
as in Subsection~\ref{subsect:H}, 
we can say that $F$ is minimal 
if $\mbox{\boldmath{$H$}}$ identically vanishes. 
\end{rem} 

Suppose that $F:M\longrightarrow E^4$ is minimal. 
Then we have 
\begin{equation*} 
F_{ww}^{\bot} =\dfrac{1}{2} (F_{uu}^{\bot} -\sqrt{-1} F_{uv}^{\bot} ), 
\end{equation*} 
where $w=u+\sqrt{-1} v$ 
and $F_{uu}^{\bot}$, $F_{uv}^{\bot}$ are the normal components 
of $F_{uu}$, $F_{uv}$ respectively. 
By Theorem~\ref{thm:GW}, $\Phi$ is represented as 
\begin{equation*} 
\begin{split} 
\Phi & =\dfrac{1}{4} ( \langle F_{uu}^{\bot} , F_{uu}^{\bot} \rangle 
                     - \langle F_{uv}^{\bot} , F_{uv}^{\bot} \rangle 
                     -2\sqrt{-1} 
                       \langle F_{uu}^{\bot} , F_{uv}^{\bot} \rangle ) \\ 
     & =(\beta^1 )^2 +(\beta^2 )^2 , 
\end{split} 
\end{equation*} 
where $\beta^1$, $\beta^2$ are as in \eqref{betagamma}. 
From Proposition~\ref{pro:b1b2}, 
we observe that $\Phi$ is holomorphic. 
Hence we obtain 

\begin{thm}\label{thm:Q} 
Let $F:M\longrightarrow E^4$ be a conformal and minimal immersion. 
Then the complex quartic differential $Q=\Phi dw^4$ 
defined as in \eqref{Qdef} is holomorphic. 
\end{thm} 

\section{Orthogonal complex structures and oriented \\ 
2-planes in $\mbox{\boldmath{$E$}}^4$}\label{sect:ocs} 

\setcounter{equation}{0} 

Let $M(4, \mbox{\boldmath{$R$}} )$ be the set of $4\times 4$ matrices 
such that the components are real numbers. 
Let $O(4)$ denote the orthogonal group of degree 4: 
\begin{equation*} 
O(4):=\{ X\in M(4, \mbox{\boldmath{$R$}} ) \ | \ {}^t\!XX=E_4\} , 
\end{equation*} 
where $E_4$ is the identity matric of degree $4$ 
and ${}^t\!X$ denotes the transposed matrix of $X$. 
Let $SO(4)$ denote the special orthogonal group of degree 4: 
\begin{equation*} 
SO(4):=\{ X\in O(4) \ | \ \det X =1\} . 
\end{equation*} 
We set 
\begin{equation} 
\Sigma :=\{ X\in O(4) \ | \ X^2 =-E_4 \} . 
\label{Sigma} 
\end{equation} 
In this survey, 
the linear transformation of $E^4$ given by an element of $\Sigma$ or 
the element itself is called 
an \textit{orthogonal complex structure\/} of $E^4$. 

In this section, 
we denote $\mbox{\boldmath{$a$}} \in E^4$ by 
\begin{equation*} 
\mbox{\boldmath{$a$}} =\left[ \begin{array}{c} 
                               a^1 \\ 
                               a^2 \\ 
                               a^3 \\ 
                               a^4 
                                \end{array} 
                       \right] . 
\end{equation*} 
Therefore, by the usual notation in linear algebra, 
for $X=[x^i_j ] \in M(4, \mbox{\boldmath{$R$}} )$, 
$X\mbox{\boldmath{$a$}}$ makes sense: 
\begin{equation*} 
X\mbox{\boldmath{$a$}} 
=\left[ \begin{array}{cccc} 
         x^1_1 & x^1_2 & x^1_3 & x^1_4 \\ 
         x^2_1 & x^2_2 & x^2_3 & x^2_4 \\ 
         x^3_1 & x^3_2 & x^3_3 & x^3_4 \\ 
         x^4_1 & x^4_2 & x^4_3 & x^4_4 
          \end{array} 
 \right] 
 \left[ \begin{array}{c} 
         a^1 \\ 
         a^2 \\ 
         a^3 \\ 
         a^4 
          \end{array} 
 \right] 
=\left[ \begin{array}{c} 
         x^1_1 a^1 +x^1_2 a^2 +x^1_3 a^3 +x^1_4 a^4 \\ 
         x^2_1 a^1 +x^2_2 a^2 +x^2_3 a^3 +x^2_4 a^4 \\ 
         x^3_1 a^1 +x^3_2 a^2 +x^3_3 a^3 +x^3_4 a^4 \\ 
         x^4_1 a^1 +x^4_2 a^2 +x^4_3 a^3 +x^4_4 a^4 
          \end{array} 
 \right] . 
\end{equation*} 
In particular, for $X\in O(4)$ 
and $\mbox{\boldmath{$a$}}$, $\mbox{\boldmath{$b$}} \in E^4$, 
\begin{equation*} 
 \langle X\mbox{\boldmath{$a$}} ,       X\mbox{\boldmath{$b$}} \rangle 
={}^t   (X\mbox{\boldmath{$a$}} )       X\mbox{\boldmath{$b$}} 
={}^t     \mbox{\boldmath{$a$}} {}^t\!X X\mbox{\boldmath{$b$}} 
={}^t     \mbox{\boldmath{$a$}} E_4      \mbox{\boldmath{$b$}} 
={}^t     \mbox{\boldmath{$a$}}          \mbox{\boldmath{$b$}} 
=\langle  \mbox{\boldmath{$a$}} ,        \mbox{\boldmath{$b$}} \rangle . 
\end{equation*} 

For $X\in \Sigma$ and $\mbox{\boldmath{$u$}} \in \mbox{\boldmath{$R$}}^4$ 
with $|\mbox{\boldmath{$u$}} |=1$, 
\begin{equation*} 
  \langle X   \mbox{\boldmath{$u$}} ,  \mbox{\boldmath{$u$}} \rangle 
= \langle X^2 \mbox{\boldmath{$u$}} , X\mbox{\boldmath{$u$}} \rangle 
=-\langle     \mbox{\boldmath{$u$}} , X\mbox{\boldmath{$u$}} \rangle 
=-\langle X   \mbox{\boldmath{$u$}} ,  \mbox{\boldmath{$u$}} \rangle , 
\end{equation*} 
and therefore $\langle X \mbox{\boldmath{$u$}} , \mbox{\boldmath{$u$}} 
               \rangle =0$. 
For $\mbox{\boldmath{$u$}}' \in E^4$ 
with $|\mbox{\boldmath{$u$}}' |=1$ 
and $\langle \mbox{\boldmath{$u$}}' ,  \mbox{\boldmath{$u$}} \rangle 
    =\langle \mbox{\boldmath{$u$}}' , X\mbox{\boldmath{$u$}} \rangle 
    =0$, 
\begin{equation*} 
\begin{split} 
&  \langle X\mbox{\boldmath{$u$}}' ,    \mbox{\boldmath{$u$}}' \rangle 
 =0, \\   
&  \langle X  \mbox{\boldmath{$u$}}' ,  \mbox{\boldmath{$u$}}  \rangle 
 = \langle X^2\mbox{\boldmath{$u$}}' , X\mbox{\boldmath{$u$}}  \rangle 
 =-\langle    \mbox{\boldmath{$u$}}' , X\mbox{\boldmath{$u$}}  \rangle 
 =0, \\ 
&  \langle X  \mbox{\boldmath{$u$}}' , X\mbox{\boldmath{$u$}}  \rangle 
 = \langle    \mbox{\boldmath{$u$}}' ,  \mbox{\boldmath{$u$}}  \rangle 
 =0. 
\end{split} 
\end{equation*}  
Therefore $\mbox{\boldmath{$u$}}$, 
         $X\mbox{\boldmath{$u$}}$, 
          $\mbox{\boldmath{$u$}}'$, 
         $X\mbox{\boldmath{$u$}}'$ form an orthonormal basis of $E^4$. 
Then we denote by $Y$ the element of $O(4)$ constructed by these vectors: 
\begin{equation*} 
Y:=[\mbox{\boldmath{$u$}}  \ 
   X\mbox{\boldmath{$u$}}  \ 
    \mbox{\boldmath{$u$}}' \ 
   X\mbox{\boldmath{$u$}}']. 
\end{equation*} 
Then the determinant $\det Y$ of $Y$ is either $1$ or $-1$. 
Now, we fix $X\in \Sigma$. 
According as $\mbox{\boldmath{$u$}}$ and $\mbox{\boldmath{$u$}}'$ 
vary so that they preserve 
\begin{equation} 
|\mbox{\boldmath{$u$}} |=|\mbox{\boldmath{$u$}}' |=1, \quad 
 \langle \mbox{\boldmath{$u$}}' ,  \mbox{\boldmath{$u$}} \rangle 
=\langle \mbox{\boldmath{$u$}}' , X\mbox{\boldmath{$u$}} \rangle 
=0, 
\label{uu'} 
\end{equation} 
$Y$ may vary. 
However, even if $\mbox{\boldmath{$u$}}$ and $\mbox{\boldmath{$u$}}'$ 
vary so that they preserve \eqref{uu'}, 
the determinant $\det Y$ does not vary, 
and therefore whether $\det Y$ is equal to $1$ or $-1$ is 
determined by $X$. 
We set 
\begin{equation*} 
\Sigma_+ :=\{ X\in \Sigma \ | \ \det Y= 1\} , \quad 
\Sigma_- :=\{ X\in \Sigma \ | \ \det Y=-1\} . 
\end{equation*} 
Then $\Sigma =\Sigma_+ \cup \Sigma_-$ and $\Sigma_+ \cap \Sigma_- =\emptyset$. 

We set 
\begin{equation*} 
I_{+, 1} 
:=\left[ \begin{array}{cccc} 
          0 & -1 & 0 &  0 \\ 
          1 &  0 & 0 &  0 \\ 
          0 &  0 & 0 & -1 \\ 
          0 &  0 & 1 &  0 
           \end{array} 
  \right] , \ 
I_{+, 2} 
:=\left[ \begin{array}{cccc} 
          0 &  0 & -1 &  0 \\ 
          0 &  0 &  0 &  1 \\ 
          1 &  0 &  0 &  0 \\ 
          0 & -1 &  0 &  0 
           \end{array} 
  \right] , \ 
I_{+, 3} 
:=\left[ \begin{array}{cccc} 
          0 &  0 &  0 & -1 \\ 
          0 &  0 & -1 &  0 \\ 
          0 &  1 &  0 &  0 \\ 
          1 &  0 &  0 &  0 
           \end{array} 
  \right]  
\end{equation*} 
and 
\begin{equation*} 
I_{-, 1} 
:=\left[ \begin{array}{cccc} 
          0 & -1 &  0 &  0 \\ 
          1 &  0 &  0 &  0 \\ 
          0 &  0 &  0 &  1 \\ 
          0 &  0 & -1 &  0 
           \end{array} 
  \right] , \ 
I_{-, 2} 
:=\left[ \begin{array}{cccc} 
          0 &  0 & -1 &  0 \\ 
          0 &  0 &  0 & -1 \\ 
          1 &  0 &  0 &  0 \\ 
          0 &  1 &  0 &  0 
           \end{array} 
  \right] , \ 
I_{-, 3} 
:=\left[ \begin{array}{cccc} 
          0 &  0 &  0 & -1 \\ 
          0 &  0 &  1 &  0 \\ 
          0 & -1 &  0 &  0 \\ 
          1 &  0 &  0 &  0 
           \end{array} 
  \right] . 
\end{equation*} 

\begin{pro}\label{pro:Iek} 
The following hold\/$:$ 
\begin{equation*} 
I_{+, 1} , I_{+, 2} , I_{+, 3} \in \Sigma_+ , \quad 
I_{-, 1} , I_{-, 2} , I_{-, 3} \in \Sigma_- . 
\end{equation*} 
\end{pro} 

We will prove 

\begin{pro}\label{pro:W+-} 
For $X\in \Sigma_{\varepsilon}$ $(\varepsilon \in \{ +, -\} )$, 
there exist $c^1$, $c^2$, $c^3 \in \mbox{\boldmath{$R$}}$ satisfying 
\begin{equation} 
X=c^1 I_{\varepsilon , 1} 
 +c^2 I_{\varepsilon , 2} 
 +c^3 I_{\varepsilon , 3} , \quad 
(c^1 )^2 +(c^2 )^2 +(c^3 )^2 =1. 
\label{c1c2c3} 
\end{equation} 
\end{pro} 

\vspace{3mm} 

\par\noindent 
\textit{Proof} \ 
For $X=[x^i_j ]\in \Sigma$, we have ${}^t\!XX=E_4$ and $X^2 =-E_4$. 
Therefore 
\begin{equation*} 
{}^t\!X ={}^t\!XE_4 =-{}^t\!XX^2 =-({}^t\!XX)X =-E_4 X =-X. 
\end{equation*} 
This means that $X$ is alternating and therefore $x^j_i =-x^i_j$: 
\begin{equation} 
X=\left[ \begin{array}{cccc} 
          0     & -x^2_1 & -x^3_1 & -x^4_1 \\ 
          x^2_1 &  0     & -x^3_2 & -x^4_2 \\ 
          x^3_1 &  x^3_2 &  0     & -x^4_3 \\ 
          x^4_1 &  x^4_2 &  x^4_3 &  0 
           \end{array} 
  \right] . 
\label{AinW} 
\end{equation} 
Noticing the first and second columns of $X$, we have 
\begin{equation} 
x^3_1 x^3_2 +x^4_1 x^4_2 =0, \quad 
(x^3_1 )^2 +(x^4_1 )^2 =(x^3_2 )^2 +(x^4_2 )^2 . 
\label{c12} 
\end{equation} 
If both sides of the second relation in \eqref{c12} are zero, 
then $X$ satisfies 
\begin{equation*} 
x^3_1 =x^3_2 =x^4_1 =x^4_2 =0, 
\end{equation*} 
and then 
$X=I_{+, 1}$ or $I_{-, 1}$. 
In the following, 
we suppose that both sides of the second relation in \eqref{c12} are 
nonzero. 
Then from \eqref{c12}, 
we see that $(x^3_1 , x^4_1 )$ is perpendicular to $(x^3_2 , x^4_2 )$ 
and that these vectors have the same length. 
Therefore we obtain 
\begin{equation} 
\left[ \begin{array}{c} 
        x^3_2 \\ 
        x^4_2 
         \end{array} 
\right] 
=
\varepsilon 
\left[ \begin{array}{c} 
        x^4_1 \\ 
       -x^3_1 
         \end{array} 
\right] 
\label{epsilon} 
\end{equation} 
for $\varepsilon =1$ or $-1$. 
In \eqref{epsilon}, suppose $\varepsilon =1$. 
Since $X\in O(4)$, we have 
\begin{equation*} 
x^2_1 x^3_1 -x^3_1 x^4_3 =0, \quad 
x^2_1 x^4_1 -x^4_1 x^4_3 =0. 
\end{equation*} 
Since at least one of $x^3_1$ and $x^4_1$ is nonzero, 
we have $x^2_1 =x^4_3$. 
Hence from \eqref{AinW}, 
we observe that $X$ is represented as 
\begin{equation*} 
X=\left[ \begin{array}{cccc} 
          0     & -x^2_1 & -x^3_1 & -x^4_1 \\ 
          x^2_1 &  0     & -x^4_1 &  x^3_1 \\ 
          x^3_1 &  x^4_1 &  0     & -x^2_1 \\ 
          x^4_1 & -x^3_1 &  x^2_1 &  0 
           \end{array} 
  \right] 
 =c^1 I_{+, 1} +c^2 I_{+, 2} +c^3 I_{+, 3} , 
\end{equation*} 
where $c^1 :=x^2_1$, $c^2 :=x^3_1$, $c^3 :=x^4_1$, 
and then $X\in \Sigma_+$. 
Similarly, if we suppose $\varepsilon =-1$, 
then $X\in \Sigma_-$. 
Thus we have proved Proposition~\ref{pro:W+-}. 
\hfill 
$\square$ 

\vspace{3mm} 

\begin{rem} 
From Proposition~\ref{pro:W+-}, 
we see that $\Sigma$ is contained in $SO(4)$. 
As is seen in the proof of Proposition~\ref{pro:W+-}, 
any element of $\Sigma$ is alternating. 
\end{rem} 

In addition to Proposition~\ref{pro:W+-}, the following holds: 

\begin{pro}\label{pro:W+-2} 
If $X\in M(4, \mbox{\boldmath{$R$}} )$ is represented as 
in \eqref{c1c2c3} for $c^1$, $c^2$, $c^3 \in \mbox{\boldmath{$R$}}$ 
satisfying $(c^1 )^2 +(c^2 )^2 +(c^3 )^2 =1$, 
then $X\in \Sigma_{\varepsilon}$. 
\end{pro} 

From Propositions~\ref{pro:W+-} and \ref{pro:W+-2}, 
we obtain 

\begin{thm}\label{thm:W+-} 
For $\varepsilon \in \{ +, -\}$, 
$\Sigma_{\varepsilon}$ is represented as 
\begin{equation*} 
\Sigma_{\varepsilon} 
  =\{ c^1 I_{\varepsilon , 1} 
    + c^2 I_{\varepsilon , 2} 
    + c^3 I_{\varepsilon , 3} \ | \ 
     (c^1 )^2 +(c^2 )^2 +(c^3 )^2 =1\} . 
\end{equation*} 
\end{thm} 

\begin{rem} 
If $X\in SO(4)$ is alternating, 
then $X$ is represented as in \eqref{c1c2c3} 
for $c^1$, $c^2$, $c^3 \in \mbox{\boldmath{$R$}}$ 
satisfying $(c^1 )^2 +(c^2 )^2 +(c^3 )^2 =1$ 
for $\varepsilon =+$ or $-$. 
Therefore $\Sigma$ coincides with the set of alternating matrices 
in $SO(4)$. 
\end{rem} 

Let $P$ be a $2$-dimensional subspace of $E^4$, 
Let $\mbox{\boldmath{$a$}}$, $\mbox{\boldmath{$b$}}$ be 
vectors of $P$ satisfying 
\begin{equation} 
|\mbox{\boldmath{$a$}} |=|\mbox{\boldmath{$b$}} |=1, \quad 
\langle \mbox{\boldmath{$a$}} , \mbox{\boldmath{$b$}} \rangle =0. 
\label{ab}
\end{equation} 
Then $\mbox{\boldmath{$a$}}$, $\mbox{\boldmath{$b$}}$ form 
an orthonormal basis of $P$. 
We fix the pair $(\mbox{\boldmath{$a$}} , 
                  \mbox{\boldmath{$b$}} )$ 
(notice that $\mbox{\boldmath{$a$}}$ is the first  vector and 
             $\mbox{\boldmath{$b$}}$ is the second vector, 
and we distinguish 
between $(\mbox{\boldmath{$a$}} , \mbox{\boldmath{$b$}} )$ 
and     $(\mbox{\boldmath{$b$}} , \mbox{\boldmath{$a$}} )$). 
Then any \textit{ordered\/} orthonormal basis of $P$ is given by 
\begin{equation} 
 (\mbox{\boldmath{$a$}}' , \ \mbox{\boldmath{$b$}}' ) 
=(\mbox{\boldmath{$a$}}  , \ \mbox{\boldmath{$b$}}  ) 
  \left[ \begin{array}{cc} 
          \cos \theta & -\sin \theta \\ 
          \sin \theta &  \cos \theta 
           \end{array} 
  \right] 
\label{ab+} 
\end{equation} 
or 
\begin{equation} 
 (\mbox{\boldmath{$a$}}' , \ \mbox{\boldmath{$b$}}' ) 
=(\mbox{\boldmath{$a$}}  , \ \mbox{\boldmath{$b$}}  ) 
  \left[ \begin{array}{cc} 
          \cos \theta &  \sin \theta \\ 
          \sin \theta & -\cos \theta 
           \end{array} 
  \right] 
\label{ab-} 
\end{equation} 
for $\theta \in [0, 2\pi )$. 

Let $\mathcal{B}_+$, $\mathcal{B}_-$ be 
the set of ordered orthonormal bases of $P$ 
represented as in \eqref{ab+}, \eqref{ab-}, respectively. 
An \textit{oriented 2-plane\/} in $E^4$ is a 2-dimensional subspace $P$ 
of $E^4$ equipped with just one of $\mathcal{B}_+$ and $\mathcal{B}_-$. 
Therefore each oriented 2-plane is represented as 
$(P, \mathcal{B}_{\varepsilon} )$ for $\varepsilon =+$ or $-$, 
and 
we distinguish between $(P, \mathcal{B}_+ )$ and $(P, \mathcal{B}_- )$. 
The set of oriented 2-planes in $E^4$ is 
denoted by $\Tilde{G}_{4, 2} (\mbox{\boldmath{$R$}} )$, 
and called an \textit{oriented real Grassmann manifold}. 

\begin{pro}\label{pro:GrW+-} 
For each $(P, \mathcal{B}_{\varepsilon} ) 
              \in \Tilde{G}_{4, 2} (\mbox{\boldmath{$R$}} )$, 
there exist unique $X_+ \in \Sigma_+$, $X_- \in \Sigma_-$ satisfying 
\begin{itemize} 
\item{$X_+ \mbox{\boldmath{$u$}} =X_- \mbox{\boldmath{$u$}} \in P$ 
\textit{for $\mbox{\boldmath{$u$}} \in P$ with\/} 
$|\mbox{\boldmath{$u$}} |=1$,} 
\item{$(\mbox{\boldmath{$u$}} , X_+ \mbox{\boldmath{$u$}} ) 
        \in \mathcal{B}_{\varepsilon}$.} 
\end{itemize} 
\end{pro} 

Proposition~\ref{pro:GrW+-} gives 
a map from $\Tilde{G}_{4, 2} (\mbox{\boldmath{$R$}} )$ 
to $\Sigma_+ \times \Sigma_-$ which assigns to 
each $(P, \mathcal{B}_{\varepsilon} ) 
          \in \Tilde{G}_{4, 2} (\mbox{\boldmath{$R$}} )$ 
a pair $(X_+ , X_- )\in \Sigma_+ \times \Sigma_-$. 

The following holds: 

\begin{pro}\label{pro:c1c2c3} 
Let $\mbox{\boldmath{$a$}} ={}^t [a^1 \ a^2 \ a^3 \ a^4 ]$, 
    $\mbox{\boldmath{$b$}} ={}^t [b^1 \ b^2 \ b^3 \ b^4 ]$ be 
vectors of $E^4$ satisfying \eqref{ab}. 
Then for each $\varepsilon \in \{ +, -\}$, 
a unique element $X_{\varepsilon}$ of $\Sigma_{\varepsilon}$ 
satisfying $X_{\varepsilon} \mbox{\boldmath{$a$}} =\mbox{\boldmath{$b$}}$ 
is given by $X_{\varepsilon} =c^1_{\varepsilon} I_{\varepsilon , 1} 
                             +c^2_{\varepsilon} I_{\varepsilon , 2} 
                             +c^3_{\varepsilon} I_{\varepsilon , 3}$, 
where 
\begin{equation*} 
\begin{split} 
c^1_{\varepsilon} & :=a^1 b^2 -a^2 b^1 +\varepsilon (a^3 b^4 -a^4 b^3 ), \\ 
c^2_{\varepsilon} & :=a^1 b^3 -a^3 b^1 +\varepsilon (a^4 b^2 -a^2 b^4 ), \\ 
c^3_{\varepsilon} & :=a^1 b^4 -a^4 b^1 +\varepsilon (a^2 b^3 -a^3 b^2 ), 
\end{split} 
\end{equation*} 
and ``$++$'', ``$+-$'' mean ``$+$'', ``$-$'' respectively. 
Moreover, $(c^1_{\varepsilon} )^2 
          +(c^2_{\varepsilon} )^2 
          +(c^3_{\varepsilon} )^2 =1$. 
\end{pro} 

Let $\mbox{\boldmath{$a$}} \wedge \mbox{\boldmath{$b$}}$ be 
an element of $M(4, \mbox{\boldmath{$R$}} )$ defined by 
\begin{equation*} 
  \mbox{\boldmath{$a$}} \wedge \mbox{\boldmath{$b$}} 
:=\mbox{\boldmath{$b$}}\,{}^t  \mbox{\boldmath{$a$}} 
 -\mbox{\boldmath{$a$}}\,{}^t  \mbox{\boldmath{$b$}} .   
\end{equation*} 
Then $\mbox{\boldmath{$b$}} \wedge \mbox{\boldmath{$a$}} 
    =-\mbox{\boldmath{$a$}} \wedge \mbox{\boldmath{$b$}}$, 
and in particular, 
$\mbox{\boldmath{$a$}} \wedge \mbox{\boldmath{$a$}}$ is the zero matrix. 
We set 
\begin{equation} 
  \mbox{\boldmath{$e$}}_1 
:=\left[ \begin{array}{c} 
   1 \\ 
   0 \\ 
   0 \\ 
   0 
           \end{array} 
  \right] , \ 
  \mbox{\boldmath{$e$}}_2 
:=\left[ \begin{array}{c} 
   0 \\ 
   1 \\ 
   0 \\ 
   0 
           \end{array} 
  \right] , \ 
  \mbox{\boldmath{$e$}}_3 
:=\left[ \begin{array}{c} 
   0 \\ 
   0 \\ 
   1 \\ 
   0 
           \end{array} 
  \right] , \ 
  \mbox{\boldmath{$e$}}_4 
:=\left[ \begin{array}{c} 
   0 \\ 
   0 \\ 
   0 \\ 
   1 
           \end{array} 
  \right] . 
\label{e1234} 
\end{equation} 
Then we obtain 
\begin{equation*} 
\begin{split} 
I_{\varepsilon , 1}
& =\mbox{\boldmath{$e$}}_1 \wedge \mbox{\boldmath{$e$}}_2 
  +\varepsilon 
   \mbox{\boldmath{$e$}}_3 \wedge \mbox{\boldmath{$e$}}_4 , \\ 
I_{\varepsilon , 2}
& =\mbox{\boldmath{$e$}}_1 \wedge \mbox{\boldmath{$e$}}_3 
  +\varepsilon 
   \mbox{\boldmath{$e$}}_4 \wedge \mbox{\boldmath{$e$}}_2 , \\ 
I_{\varepsilon , 3}
& =\mbox{\boldmath{$e$}}_1 \wedge \mbox{\boldmath{$e$}}_4 
  +\varepsilon 
   \mbox{\boldmath{$e$}}_2 \wedge \mbox{\boldmath{$e$}}_3 . 
\end{split} 
\end{equation*} 
We set 
\begin{equation*} 
  \textstyle\bigwedge^2\!E^4 
:=\{ \mbox{\boldmath{$a$}} \wedge \mbox{\boldmath{$b$}} \ | \ 
     \mbox{\boldmath{$a$}} ,      \mbox{\boldmath{$b$}} 
     \in E^4 \} . 
\end{equation*} 
Then $\bigwedge^2\!E^4$ is a $6$-dimensional subspace of 
a $16$-dimensional vector space $M(4, \mbox{\boldmath{$R$}} )$ 
spanned by $I_{\pm , 1}$, $I_{\pm , 2}$, $I_{\pm , 3}$, 
and $\bigwedge^2\!E^4$ coincides with the set of alternating matrices 
in $M(4, \mbox{\boldmath{$R$}} )$. 
For $\varepsilon \in \{ +, -\}$, 
let $\bigwedge^2_{\varepsilon}\!E^4$ be the subspaces of $\bigwedge^2\!E^4$ 
spanned by $I_{\varepsilon , 1}$, 
           $I_{\varepsilon , 2}$, 
           $I_{\varepsilon , 3}$. 
Then in terms of the standard inner product 
of $M(4, \mbox{\boldmath{$R$}} )$ given by 
\begin{equation*} 
\langle X, Y\rangle :=\dfrac{1}{4} {\rm tr}\,({}^t\!XY) \quad 
(X, Y\in M(4, \mbox{\boldmath{$R$}} )), 
\end{equation*} 
$\bigwedge^2_+\!E^4$ is perpendicular to $\bigwedge^2_-\!E^4$, 
and $\Sigma_{\varepsilon}$ is considered to be the unit sphere 
in $\bigwedge^2_{\varepsilon}\!E^4$ centered at the origin. 

\begin{rem} 
We have $\Sigma_{\varepsilon} 
        =\bigwedge^2_{\varepsilon}\!E^4 \cap SO(4)$. 
\end{rem} 

We set 
\begin{equation*} 
H_1 :=\left\{ \left[ \begin{array}{cccc} 
                      b_1 & -b_2 & -b_3 & -b_4 \\ 
                      b_2 &  b_1 &  b_4 & -b_3 \\ 
                      b_3 & -b_4 &  b_1 &  b_2 \\ 
                      b_4 &  b_3 & -b_2 &  b_1 
                       \end{array} 
              \right] \in SO(4)\right\} . 
\end{equation*} 
Then $H_1$ is isomorphic to the special unitary group of degree $2$: 
\begin{equation*} 
SU(2):=\left\{ \left. 
               \left[ \begin{array}{cc} 
                       \alpha & -\overline{\beta} \\ 
                       \beta  &  \overline{\alpha} 
                        \end{array} 
               \right] \ \right| \ 
               \alpha , \beta \in \mbox{\boldmath{$C$}}, \ 
              |\alpha |^2 +|\beta |^2 =1 
       \right\} . 
\end{equation*} 
We set 
\begin{equation*} 
H_2 :=\left\{ \left[ \begin{array}{cccc} 
                      1 & 0      & 0      & 0      \\ 
                      0 & c_{11} & c_{12} & c_{13} \\ 
                      0 & c_{21} & c_{22} & c_{23} \\ 
                      0 & c_{31} & c_{32} & c_{33} 
                       \end{array} 
              \right] \in SO(4)\right\} . 
\end{equation*} 
Then for each element $X$ of $SO(4)$, 
there exist unique elements $B\in H_1$, $C\in H_2$ satisfying $X=BC$. 
If we represent $B\in H_1$ 
as $B=[\mbox{\boldmath{$b$}}_1 \ 
       \mbox{\boldmath{$b$}}_2 \ 
       \mbox{\boldmath{$b$}}_3 \ 
       \mbox{\boldmath{$b$}}_4 ]$, 
then we obtain 
\begin{equation*} 
\begin{split} 
    \mbox{\boldmath{$b$}}_1 \wedge \mbox{\boldmath{$b$}}_2 
   +\mbox{\boldmath{$b$}}_3 \wedge \mbox{\boldmath{$b$}}_4 
& = I_{+, 1} , \\ 
    \mbox{\boldmath{$b$}}_1 \wedge \mbox{\boldmath{$b$}}_3 
   +\mbox{\boldmath{$b$}}_4 \wedge \mbox{\boldmath{$b$}}_2 
& = I_{+, 2} , \\ 
    \mbox{\boldmath{$b$}}_1 \wedge \mbox{\boldmath{$b$}}_4 
   +\mbox{\boldmath{$b$}}_2 \wedge \mbox{\boldmath{$b$}}_3 
& = I_{+, 3} 
\end{split} 
\end{equation*} 
and 
\begin{equation*} 
\begin{split} 
& (\mbox{\boldmath{$b$}}_1 \wedge \mbox{\boldmath{$b$}}_2 
  -\mbox{\boldmath{$b$}}_3 \wedge \mbox{\boldmath{$b$}}_4 , \ 
   \mbox{\boldmath{$b$}}_1 \wedge \mbox{\boldmath{$b$}}_3 
  -\mbox{\boldmath{$b$}}_4 \wedge \mbox{\boldmath{$b$}}_2 , \ 
   \mbox{\boldmath{$b$}}_1 \wedge \mbox{\boldmath{$b$}}_4 
  -\mbox{\boldmath{$b$}}_2 \wedge \mbox{\boldmath{$b$}}_3 ) \\ 
& =(I_{-, 1} , \ I_{-, 2} , \ I_{-, 3} )\Tilde{B} , 
\end{split} 
\end{equation*} 
where 
\begin{equation*} 
\Tilde{B} 
:=\left[ \begin{array}{ccc} 
          b^2_1 +b^2_2 - b^2_3 -b^2_4 &  2b_1    b_4   +2b_2    b_3 
                                      & -2b_1    b_3   +2b_2    b_4 \\ 
        -2b_1    b_4   +2b_2    b_3   &   b^2_1 +b^2_3 - b^2_2 -b^2_4 
                                      &  2b_1    b_2   +2b_3    b_4 \\  
         2b_1    b_3   +2b_2    b_4   & -2b_1    b_2   +2b_3    b_4 
                                      &   b^2_1 +b^2_4 - b^2_2 -b^2_3 
           \end{array} 
  \right] . 
\end{equation*} 
The matrix $\Tilde{B}$ is an element of $SO(3)$, 
and therefore we obtain a map $\Phi$ from $H_1$ to $SO(3)$ 
by $\Phi (B):=\Tilde{B}$. 
We see that $\Phi :H_1 \longrightarrow SO(3)$ is a homomorphism. 
In addition, $\Phi$ is a double covering from $H_1$ onto $SO(3)$. 
If we represent $C\in H_2$ 
as $C=[\mbox{\boldmath{$c$}}_1 \ 
       \mbox{\boldmath{$c$}}_2 \ 
       \mbox{\boldmath{$c$}}_3 \ 
       \mbox{\boldmath{$c$}}_4 ]$, 
then we obtain 
\begin{equation*} 
\begin{split} 
& (\mbox{\boldmath{$c$}}_1 \wedge \mbox{\boldmath{$c$}}_2 
  +\varepsilon 
   \mbox{\boldmath{$c$}}_3 \wedge \mbox{\boldmath{$c$}}_4 , \ 
   \mbox{\boldmath{$c$}}_1 \wedge \mbox{\boldmath{$c$}}_3 
  +\varepsilon 
   \mbox{\boldmath{$c$}}_4 \wedge \mbox{\boldmath{$c$}}_2 , \ 
   \mbox{\boldmath{$c$}}_1 \wedge \mbox{\boldmath{$c$}}_4 
  +\varepsilon 
   \mbox{\boldmath{$c$}}_2 \wedge \mbox{\boldmath{$c$}}_3 ) \\ 
& =(I_{\varepsilon , 1} , \ I_{\varepsilon , 2} , \ I_{\varepsilon , 3} ) 
   \left[ \begin{array}{ccc} 
           c_{11} & c_{12} & c_{13} \\ 
           c_{21} & c_{22} & c_{23} \\ 
           c_{31} & c_{32} & c_{33} 
            \end{array} 
   \right] . 
\end{split} 
\end{equation*} 
Therefore we obtain 
a double covering homomorphism $\Tilde{\Phi} : SO(4)\longrightarrow 
                                               SO(3)\times SO(3)$ 
by $\Tilde{\Phi} (X):=(C, \Tilde{B} C)$, 
and by $\Tilde{\Phi}$, 
we observe that for each pair $(X_+ , X_- )\in \Sigma_+ \times \Sigma_-$, 
there exists a unique oriented 2-plane $(P, \mathcal{B}_{\varepsilon} )$ 
which corresponds to $(X_+ , X_- )$ by Proposition~\ref{pro:GrW+-}. 
Therefore we obtain 

\begin{thm}\label{thm:GrW+-} 
The map given in Proposition~\ref{pro:GrW+-} yields 
a one-to-one correspondence 
between $\Tilde{G}_{4, 2} (\mbox{\boldmath{$R$}} )$ 
and     $\Sigma_+ \times \Sigma_-$. 
\end{thm} 

\section{The twistor lifts and the Gauss maps of surfaces} 

\setcounter{equation}{0} 

Let $F:M\longrightarrow E^4$ be a conformal immersion. 
Then, as in Subsection~\ref{subsect:H}, 
on each holomorphic coordinate neighborhood $(U, w)$ of $M$ 
with $w=u+\sqrt{-1} v$, 
we have $g_{12} =0$, and 
there exists a $C^{\infty}$-function $\alpha$ 
on $U$ satisfying $e^{2\alpha} =g_{11} =g_{22}$. 
Since $g_{ij}$ ($i, j=1, 2$) satisfy \eqref{gij}, 
\begin{equation} 
\mbox{\boldmath{$t$}}_1 :=\dfrac{1}{e^{\alpha}} F_u , \quad 
\mbox{\boldmath{$t$}}_2 :=\dfrac{1}{e^{\alpha}} F_v 
\label{t1t2} 
\end{equation} 
give an orthonormal basis of $T_p M$ for each point $p$ of $U$. 
Suppose that for $P:=dF (T_p M)$, 
$\mathcal{B}_+$ contains $(\mbox{\boldmath{$t$}}_1 , 
                           \mbox{\boldmath{$t$}}_2 )$. 
We can identify $dF(T_p M)$ with a $2$-dimensional subspace 
of $E^4$ by natural parallelism. 
Therefore we obtain a map $G$ from $U$ 
into $\Tilde{G}_{4, 2} (\mbox{\boldmath{$R$}} )$ 
by $G(p):=(dF(T_p M), \mathcal{B}_+ )$ for $p\in U$. 
In addition, since $M$ is orientable, 
$G$ is well-defined on $M$. 
By Proposition~\ref{pro:GrW+-}, 
$G(p)=(dF(T_p M), \mathcal{B}_+ )$ gives 
a pair $(F_+ (p), F_- (p))\in \Sigma_+ \times \Sigma_-$. 
We call 
\begin{center} 
$G           :M\longrightarrow \Tilde{G}_{4, 2} (\mbox{\boldmath{$R$}} )$ \ 
or \ 
$(F_+ , F_- ):M\longrightarrow \Sigma_+ \times \Sigma_-$ 
\end{center} 
the \textit{Gauss map\/} of $F$, and 
we call each of $F_{\pm} :M\longrightarrow \Sigma_{\pm}$ 
a \textit{twistor lift\/} of $F$. 

On $(U, w)$, we set 
\begin{equation*} 
\psi^i :=\dfrac{\partial f^i}{\partial w} 
        =\dfrac{1}{2} 
         \left(    \dfrac{\partial f^i}{\partial u} 
        -\sqrt{-1} \dfrac{\partial f^i}{\partial v} 
         \right) . 
\end{equation*} 
Then 
\begin{itemize} 
\item{if $F$ is minimal, then each $\psi^i$ is holomorphic, 
since $f^i$ is harmonic,} 
\item{$\sum^4_{i=1} |\psi^i |^2 =e^{2\alpha} /2$,} 
\item{$\sum^4_{i=1} (\psi^i )^2 =0$, 
since $(u, v)$ are isothermal coordinates.} 
\end{itemize} 
Referring to Proposition~\ref{pro:c1c2c3}, 
we obtain 

\begin{pro}\label{pro:tl} 
On $(U, w)$, 
the twistor lifts $F_{\pm}$ are represented as 
\begin{equation*} 
F_{\varepsilon} =-\dfrac{2\sqrt{-1}}{e^{2\alpha}} 
                  \displaystyle\sum^3_{k=1} \Psi^k_{\varepsilon} 
                                                 I_{\varepsilon , k} \ \ 
(\varepsilon =+, -), 
\end{equation*} 
where 
\begin{equation*} 
\begin{split} 
\Psi^1_{\varepsilon} 
& :=\psi^1 \overline{\psi}^2 -\psi^2 \overline{\psi}^1 
   +\varepsilon 
   (\psi^3 \overline{\psi}^4 -\psi^4 \overline{\psi}^3 ), \\ 
\Psi^2_{\varepsilon} 
& :=\psi^1 \overline{\psi}^3 -\psi^3 \overline{\psi}^1 
   +\varepsilon  
   (\psi^4 \overline{\psi}^2 -\psi^2 \overline{\psi}^4 ), \\ 
\Psi^3_{\varepsilon} 
& :=\psi^1 \overline{\psi}^4 -\psi^4 \overline{\psi}^1 
   +\varepsilon 
   (\psi^2 \overline{\psi}^3 -\psi^3 \overline{\psi}^2 )  \\ 
\end{split} 
\end{equation*} 
and $\overline{\psi}^i$ is the conjugate of $\psi^i :$ 
\begin{equation*} 
\overline{\psi}^i 
:=\dfrac{\partial f^i}{\partial \overline{w}} 
 =\dfrac{1}{2} 
  \left(    \dfrac{\partial f^i}{\partial u} 
 +\sqrt{-1} \dfrac{\partial f^i}{\partial v} 
  \right) . 
\end{equation*} 
\end{pro} 

\begin{rem} 
Functions $-(2\sqrt{-1} /e^{2\alpha} )\Psi^k_{\varepsilon}$ 
($k=1, 2, 3$) are real-valued and satisfy 
\begin{equation} 
\displaystyle\sum^3_{k=1} 
\left( -\dfrac{2\sqrt{-1}}{e^{2\alpha}} \Psi^k_{\varepsilon} 
\right)^2 =1. 
\label{sum2} 
\end{equation} 
\end{rem} 

Suppose $-(2\sqrt{-1} /e^{2\alpha} )\Psi^3_{\varepsilon} \not= 1$ on $U$. 
Then by \eqref{sum2}, 
there exist complex-valued functions $g_{\pm}$ on $U$ satisfying 
\begin{equation} 
   \dfrac{1}{|g_{\varepsilon} |^2 +1} 
 (2\,{\rm Re}\,g_{\varepsilon} , \ 
  2\,{\rm Im}\,g_{\varepsilon} , \ 
              |g_{\varepsilon} |^2 -1)  
 =-\dfrac{2\sqrt{-1}}{e^{2\alpha}} 
   (\Psi^1_{\varepsilon} , \Psi^2_{\varepsilon} , \Psi^3_{\varepsilon} ) 
\label{gsp} 
\end{equation} 
for $\varepsilon =+, -$. 

\begin{rem} 
Let $\pi_+$, $p_{\pm}$ be as in Example~\ref{ex:Rsphere}. 
We set $(a, b):=\pi_+ (x, y, z)$ 
for $(x, y, z)\in U_+$, and $c:=a+\sqrt{-1} b$. 
Then by \eqref{pipm}, we have 
\begin{equation*} 
(x, y, z)=\dfrac{1}{|c|^2 +1} (2\,{\rm Re}\,c, \ 
                               2\,{\rm Im}\,c, \ 
                                           |c|^2 -1). 
\end{equation*} 
Therefore the stereographic projection $\pi_+$ maps  
the point of $U_+$ given by the right hand side 
of \eqref{gsp} 
to the point $({\rm Re}\,g_{\varepsilon} , {\rm Im}\,g_{\varepsilon} )$ 
of $\mbox{\boldmath{$R$}}^2$. 
\end{rem} 

The following is a result given in \cite{HO}. 

\begin{thm}\label{thm:HO} 
If $F$ is minimal, 
then the complex-valued functions $g_+$, $\overline{g}_-$ are holomorphic. 
\end{thm} 

\vspace{3mm} 

\par\noindent 
\textit{Outline of the proof} \ 
This outline is based on \cite{ando(7)}. 
We set $\varepsilon =+$. 
From $\sum^4_{i=1} |\psi^i |^2 =e^{2\alpha} /2$ 
and \eqref{gsp}, we obtain 
\begin{equation} 
-\dfrac{\sqrt{-1}}{\sum^4_{i=1} |\psi^i |^2} 
(\Psi^1_+ +\sqrt{-1} \Psi^2_+ ) 
=\left( 1+\dfrac{\sqrt{-1}}{\sum^4_{i=1} |\psi^i |^2} \Psi^3_+ 
 \right) g_+ . 
\label{gsp2} 
\end{equation} 
Then \eqref{gsp2} is rewritten into 
\begin{equation} 
\begin{split} 
&  (          \psi^1  +\sqrt{-1}            \psi^4  ) 
   (\overline{\psi}^2 +\sqrt{-1}\,\overline{\psi}^3 ) 
  -(          \psi^2  +\sqrt{-1}            \psi^3  ) 
   (\overline{\psi}^1 +\sqrt{-1}\,\overline{\psi}^4 ) \\ 
& = \sqrt{-1} g_+ 
  ((          \psi^1  -\sqrt{-1}            \psi^4  ) 
   (\overline{\psi}^1 +\sqrt{-1}\,\overline{\psi}^4 ) \\ 
&   \quad \quad \quad \quad \quad \quad 
  +(          \psi^2  -\sqrt{-1}            \psi^3  ) 
   (\overline{\psi}^2 +\sqrt{-1}\,\overline{\psi}^3 )).  
\end{split} 
\label{gsp3} 
\end{equation} 
Notice 
\begin{equation*} 
 (          \psi^1  -\sqrt{-1}            \psi^4  ) 
 (\overline{\psi}^1 +\sqrt{-1}\,\overline{\psi}^4 ) 
+(          \psi^2  -\sqrt{-1}            \psi^3  ) 
 (\overline{\psi}^2 +\sqrt{-1}\,\overline{\psi}^3 ) 
  \not= 0, 
\end{equation*} 
because of $-(2\sqrt{-1} /e^{2\alpha} )\Psi^3_+ \not= 1$. 
Since $\sum^4_{i=1} (\psi^i )^2 =0$, we have 
\begin{equation} 
  (\psi^2 +\sqrt{-1} \psi^3 )(\psi^2 -\sqrt{-1} \psi^3 ) 
=-(\psi^1 +\sqrt{-1} \psi^4 )(\psi^1 -\sqrt{-1} \psi^4 ).  
\label{psi2314} 
\end{equation} 
Applying \eqref{psi2314} to \eqref{gsp3}, we obtain 
\begin{equation*} 
g_+ =\dfrac{1}{\sqrt{-1}}\,
     \dfrac{\psi^1 +\sqrt{-1} \psi^4}{\psi^2 -\sqrt{-1} \psi^3} \ \ 
{\rm or} \ \ 
     \sqrt{-1}\,
     \dfrac{\psi^2 +\sqrt{-1} \psi^3}{\psi^1 -\sqrt{-1} \psi^4} . 
\end{equation*} 
Since $\psi^i$ ($i=1, 2, 3, 4$) are holomorphic, 
we see that $g_+$ is holomorphic. 
In the case of $\varepsilon =-$, 
we similarly obtain the result for $\overline{g}_-$. 
\hfill 
$\square$ 

\vspace{3mm} 

\begin{rem} 
In the above discussion, 
we supposed $-(2\sqrt{-1} /e^{2\alpha} )\Psi^3_{\varepsilon} \not= 1$. 
In the case where we don't suppose this condition, 
Theorem~\ref{thm:HO} can be extended to 
the result that $g_+$, $\overline{g}_-$ are holomorphic maps 
from $M$ into $\mbox{\boldmath{$C$}}\!P^1$. 
\end{rem} 

The following is a result based on \cite{friedrich} 
in the case where the ambient space is $E^4$. 

\begin{thm}\label{thm:F} 
Suppose that $F$ is minimal. 
Then the following are mutually equivalent\/$:$ 
\begin{itemize} 
\item[{\rm (a)}]{the holomorphic quartic differential $Q$ 
as in \eqref{Qdef} is identically zero\/$;$} 
\item[{\rm (b)}]{$(b^1_{11} )^2 -(b^1_{12} )^2 
                 +(b^2_{11} )^2 -(b^2_{12} )^2 =0$ 
and $b^1_{11} b^1_{12} +b^2_{11} b^2_{12} =0$, 
where $b^k_{ij} (i, j, k=1, 2)$ are as in \eqref{Fperpab}\/$;$} 
\item[{\rm (c)}]{$(b^1_{11} )^2 +(b^1_{12} )^2 
                 -(b^2_{11} )^2 -(b^2_{12} )^2 =0$ 
and $b^1_{11} b^2_{11} +b^1_{12} b^2_{12} =0;$} 
\item[{\rm (d)}]{at each point of $M$, 
the eigenvalues of $(\cos \theta )A_1 +(\sin \theta )A_2$ 
with \eqref{ntopka2} do not depend on 
the choice of $\theta \in \mbox{\boldmath{$R$}} ;$} 
\item[{\rm (e)}]{one of the twistor lifts $F_{\pm}$ is constant.} 
\end{itemize} 
\end{thm} 

A conformal and minimal immersion $F$ is said to be \textit{isotropic\/} 
if $F$ satisfies one of (a)$\sim$(e) in Theorem~\ref{thm:F}. 

\vspace{3mm} 

\par\noindent 
\textit{Proof of Theorem~\ref{thm:F}} \ 
We immediately see that (a) is equivalent to (b). 
In addition, we also see that (b) is equivalent to (c). 
By Propositions~\ref{pro:A} and \ref{pro:ic}, we obtain 
\begin{equation} 
\begin{split} 
& e^{2\alpha} ((\cos \theta )A_1 +(\sin \theta )A_2 ) \\ 
& =\cos \theta \left[ 
               \begin{array}{cc} 
                b^1_{11} &  b^1_{12} \\ 
                b^1_{21} & -b^1_{11} 
                 \end{array} 
               \right]  
  +\sin \theta \left[ 
               \begin{array}{cc} 
                b^2_{11} &  b^2_{12} \\ 
                b^2_{21} & -b^2_{11} 
                 \end{array} 
               \right] . 
\end{split} 
\label{thetaA} 
\end{equation} 
Since the trace of the matrix in the right hand side of \eqref{thetaA} is 
zero, 
the eigenvalues of the matrix at each point of $M$ are 
given by $\pm \lambda$ 
for a real number $\lambda$. 
Since the matrix is symmetric, we have 
\begin{equation*} 
\begin{split} 
\lambda^2 = &  ((\cos \theta )b^1_{11} +(\sin \theta )b^2_{11} )^2 
              +((\cos \theta )b^1_{12} +(\sin \theta )b^2_{12} )^2 \\ 
          = &   ((b^1_{11} )^2 +(b^1_{12} )^2 )\cos^2 \theta  
              +2(b^1_{11} b^2_{11} +b^1_{12} b^2_{12} ) 
                 \cos \theta \sin \theta \\ 
            & + ((b^2_{11} )^2 +(b^2_{12} )^2 )\sin^2 \theta \\ 
          = &  \dfrac{1}{2} ((b^1_{11} )^2 +(b^1_{12} )^2 
                           + (b^2_{11} )^2 +(b^2_{12} )^2 ) \\ 
            & +\dfrac{1}{2} ((b^1_{11} )^2 +(b^1_{12} )^2 
                           - (b^2_{11} )^2 -(b^2_{12} )^2 )\cos 2\theta \\ 
            & +(b^1_{11} b^2_{11} +b^1_{12} b^2_{12} )\sin 2\theta .  
\end{split} 
\end{equation*} 
Therefore $\lambda$ does not depend on the choice of $\theta$ 
if and only if (c) holds. 
Therefore (d) is equivalent to (c). \\ 
Let $\mbox{\boldmath{$n$}}_1$, $\mbox{\boldmath{$n$}}_2$ be 
$\mbox{\boldmath{$R$}}^4$-valued $C^{\infty}$-functions on $U$ 
satisfying \eqref{n1n2}. 
Then $\mbox{\boldmath{$t$}}_1$, $\mbox{\boldmath{$t$}}_2$, 
     $\mbox{\boldmath{$n$}}_1$, $\mbox{\boldmath{$n$}}_2$ form  
an orthonormal basis of $\mbox{\boldmath{$R$}}^4$ at each point of $U$ 
and 
we can suppose that the element of $O(4)$ given by 
$[\mbox{\boldmath{$t$}}_1 \ \mbox{\boldmath{$t$}}_2 \ 
  \mbox{\boldmath{$n$}}_1 \ \mbox{\boldmath{$n$}}_2 ]$ 
is contained in $SO(4)$, that is, 
$\det\,[\mbox{\boldmath{$t$}}_1 \ \mbox{\boldmath{$t$}}_2 \ 
        \mbox{\boldmath{$n$}}_1 \ \mbox{\boldmath{$n$}}_2 ]=1$. 
Then a twistor lift $F_{\varepsilon}$ of $F$ is represented as 
\begin{equation*} 
 F_{\varepsilon} 
=\mbox{\boldmath{$t$}}_1 \wedge \mbox{\boldmath{$t$}}_2 
+\varepsilon 
 \mbox{\boldmath{$n$}}_1 \wedge \mbox{\boldmath{$n$}}_2 . 
\end{equation*} 
By Theorem~\ref{thm:GW} and \eqref{t1t2}, we obtain 
\begin{equation*} 
\begin{split} 
F_{+u} = &  
 \mbox{\boldmath{$t$}}_{1u} \wedge \mbox{\boldmath{$t$}}_2 
+\mbox{\boldmath{$t$}}_1    \wedge \mbox{\boldmath{$t$}}_{2u} 
+\mbox{\boldmath{$n$}}_{1u} \wedge \mbox{\boldmath{$n$}}_2 
+\mbox{\boldmath{$n$}}_1    \wedge \mbox{\boldmath{$n$}}_{2u} \\ 
       = &  
 (-\alpha_v \mbox{\boldmath{$t$}}_2 
  + e^{-\alpha} ( b^1_{11} \mbox{\boldmath{$n$}}_1 
                 +b^2_{11} \mbox{\boldmath{$n$}}_2 ))\wedge 
 \mbox{\boldmath{$t$}}_2 \\ 
         & 
+\mbox{\boldmath{$t$}}_1    \wedge 
 ( \alpha_v \mbox{\boldmath{$t$}}_1 
  + e^{-\alpha} ( b^1_{12} \mbox{\boldmath{$n$}}_1 
                 +b^2_{12} \mbox{\boldmath{$n$}}_2 )) \\ 
         & 
+(- e^{-\alpha} ( b^1_{11} \mbox{\boldmath{$t$}}_1 
                 +b^1_{12} \mbox{\boldmath{$t$}}_2 ) 
  +\gamma_1 \mbox{\boldmath{$n$}}_2 )\wedge  
 \mbox{\boldmath{$n$}}_2 \\ 
         & 
+\mbox{\boldmath{$n$}}_1    \wedge 
 (- e^{-\alpha} ( b^2_{11} \mbox{\boldmath{$t$}}_1 
                 +b^2_{12} \mbox{\boldmath{$t$}}_2 ) 
  -\gamma_1 \mbox{\boldmath{$n$}}_1 ) \\ 
       = &  e^{-\alpha} ( b^2_{11} +b^1_{12} ) 
          ( \mbox{\boldmath{$t$}}_1 \wedge \mbox{\boldmath{$n$}}_1 
           +\mbox{\boldmath{$n$}}_2 \wedge \mbox{\boldmath{$t$}}_2 ) \\ 
         & +e^{-\alpha} (-b^1_{11} +b^2_{12} ) 
          ( \mbox{\boldmath{$t$}}_1 \wedge \mbox{\boldmath{$n$}}_2 
           +\mbox{\boldmath{$t$}}_2 \wedge \mbox{\boldmath{$n$}}_1 ), 
\end{split} 
\end{equation*} 
and by similar computations, we also obtain 
\begin{equation*} 
\begin{split} 
F_{+v} = & 
 \mbox{\boldmath{$t$}}_{1v} \wedge \mbox{\boldmath{$t$}}_2 
+\mbox{\boldmath{$t$}}_1    \wedge \mbox{\boldmath{$t$}}_{2v} 
+\mbox{\boldmath{$n$}}_{1v} \wedge \mbox{\boldmath{$n$}}_2 
+\mbox{\boldmath{$n$}}_1    \wedge \mbox{\boldmath{$n$}}_{2v} \\ 
       = &  e^{-\alpha} ( b^1_{22} +b^2_{21} ) 
            ( \mbox{\boldmath{$t$}}_1 \wedge \mbox{\boldmath{$n$}}_1 
             +\mbox{\boldmath{$n$}}_2 \wedge \mbox{\boldmath{$t$}}_2 ) \\ 
         & +e^{-\alpha} (-b^1_{21} +b^2_{22} ) 
            ( \mbox{\boldmath{$t$}}_1 \wedge \mbox{\boldmath{$n$}}_2 
             +\mbox{\boldmath{$t$}}_2 \wedge \mbox{\boldmath{$n$}}_1 ). 
\end{split} 
\end{equation*} 
Therefore $F_+$ is constant if and only if 
\begin{equation} 
b^1_{11} =b^2_{12} =-b^1_{22} , \quad 
b^2_{22} =b^1_{12} =-b^2_{11} , 
\label{bij} 
\end{equation} 
and \eqref{bij} yields (c). 
Similarly, we can show that $F_-$ is constant if and only if 
\begin{equation} 
b^1_{11} =-b^2_{12} =-b^1_{22} , \quad 
b^2_{22} =-b^1_{12} =-b^2_{11} , 
\label{bij2} 
\end{equation} 
and \eqref{bij2} yields (c). 
Moreover, if $F$ is minimal so that (c) holds, 
then either \eqref{bij} or \eqref{bij2} holds, 
and therefore one of $F_{\pm}$ is constant. 
Hence (e) is equivalent to (c), 
and therefore we obtain Theorem~\ref{thm:F}. 
\hfill 
$\square$ 

\vspace{3mm} 

\begin{rem} 
Since as in Remark~\ref{rem:rmso}, 
$(\cos \theta )A_1 +(\sin \theta )A_2$ is 
the representation matrix of the shape operator 
for the normal vector field $\mbox{\boldmath{$n$}} (\theta )$ 
with length one, 
the condition (d) in Theorem~\ref{thm:F} just means that 
the principal curvatures of $F$ 
with respect to $\mbox{\boldmath{$n$}} (\theta )$ 
do not depend on the choice of $\theta \in \mbox{\boldmath{$R$}}$. 
\end{rem} 

\begin{rem} 
Let $X$ be an orthogonal complex structure of $E^4$, i.e., 
an element of $\Sigma$ (recall \eqref{Sigma}). 
A conformal immersion $F:M\longrightarrow E^4$ or its image is called 
a \textit{complex curve\/} with respect to 
the orthogonal complex structure $X$ 
if $F$ satisfies $F_v =XF_u$ on $M$. 
If $F$ is a complex curve with respect to $X\in \Sigma$, 
then by 
\begin{equation*} 
F_{vv} =XF_{uv} =X^2 F_{uu} =-F_{uu} , 
\end{equation*} 
$F$ is minimal. 
If the twistor lift $F_{\varepsilon}$ is constant 
for $\varepsilon =+$ or $-$, 
then we have $F_v =F_{\varepsilon} F_u$ on $M$ 
and therefore $F$ is a complex curve 
with respect to $X:=F_{\varepsilon}$. 
If $F$ is a complex curve with respect to $X\in \Sigma_{\varepsilon}$ 
for $\varepsilon =+$ or $-$, 
then the twistor lift $F_{\varepsilon}$ is given by $X$ 
and therefore $F_{\varepsilon}$ is constant. 
Hence the condition (e) in Theorem~\ref{thm:F} just means that 
$F$ is a complex curve with respect to 
an orthogonal complex structure of $E^4$. 
A complex curve in $E^4$ is characterized in terms of 
a relation between the induced metric and a holomorphic cubic differential 
(holomorphic $3$-differential) on $M$ (\cite{ando(1)}). 
\end{rem} 

\section{Cases where the ambient spaces are not 
\mbox{\boldmath{$E^4$}}}\label{sect:notE4} 

\setcounter{equation}{0} 

\subsection{Riemannian cases} 

In an oriented $4$-dimensional Riemannian manifold $N$, 
the definition of the twistor lifts of a surface needs 
the twistor spaces associated with $N$. 

Let $V$ be a $4$-dimensional vector space. 
Let $V\otimes V$ be the tensor product of $V$ and itself, 
which is a vector space. 
There exists 
a nondegenerate bilinear map $\Phi : V\times V\longrightarrow V\otimes V$ 
satisfying the universality of the tensor product, 
and $V\otimes V$ is generated by the image of $\Phi$. 
For $\mbox{\boldmath{$a$}}$, $\mbox{\boldmath{$b$}} \in V$, 
$\Phi (\mbox{\boldmath{$a$}} , \mbox{\boldmath{$b$}} )$ is 
denoted by $\mbox{\boldmath{$a$}} \otimes \mbox{\boldmath{$b$}}$. 
Then, for a basis $\{ \xi_1 , \xi_2 , \xi_3 , \xi_4 \}$ of $V$, 
$\xi_i \otimes \xi_j$ ($i, j=1, 2, 3, 4$) generate $V\otimes V$. 
In addition, 
by the universality of the tensor product, 
they form a basis of $V\otimes V$, and 
by the map $\xi_i \otimes \xi_j \longmapsto 
            \mbox{\boldmath{$e$}}_j {}^t \mbox{\boldmath{$e$}}_i$ 
with \eqref{e1234}, 
$V\otimes V$ is isomorphic to $M(4, \mbox{\boldmath{$R$}} )$. 
Since $\Phi$ is bilinear, we have 
\begin{equation*} 
\begin{split} 
   (  \mbox{\boldmath{$a$}} +\mbox{\boldmath{$b$}} )\otimes 
      \mbox{\boldmath{$c$}} 
& =   \mbox{\boldmath{$a$}}  \otimes   \mbox{\boldmath{$c$}} 
   +  \mbox{\boldmath{$b$}}  \otimes   \mbox{\boldmath{$c$}} , \\ 
      \mbox{\boldmath{$a$}}  \otimes 
   (  \mbox{\boldmath{$b$}} +\mbox{\boldmath{$c$}} )
& =   \mbox{\boldmath{$a$}}  \otimes   \mbox{\boldmath{$b$}} 
   +  \mbox{\boldmath{$a$}}  \otimes   \mbox{\boldmath{$c$}} , \\ 
   (r \mbox{\boldmath{$a$}}) \otimes   \mbox{\boldmath{$b$}} 
& =   \mbox{\boldmath{$a$}}  \otimes (r\mbox{\boldmath{$b$}} ) 
  = r(\mbox{\boldmath{$a$}}  \otimes   \mbox{\boldmath{$b$}} ) 
\end{split} 
\end{equation*} 
for  $\mbox{\boldmath{$a$}}$, 
     $\mbox{\boldmath{$b$}}$, 
     $\mbox{\boldmath{$c$}} \in V$, 
$r\in \mbox{\boldmath{$R$}}$. 
Let  $\bigwedge^2 V$ be the two-fold exterior power of $V$. 
Then $\bigwedge^2 V$ is a subspace of $V\otimes V$ 
generated by $\{ \mbox{\boldmath{$a$}} \wedge  \mbox{\boldmath{$b$}} \ | \ 
                 \mbox{\boldmath{$a$}} ,       \mbox{\boldmath{$b$}} \in V\}$ 
with 
\begin{equation*} 
  \mbox{\boldmath{$a$}} \wedge  \mbox{\boldmath{$b$}} 
:=\mbox{\boldmath{$a$}} \otimes \mbox{\boldmath{$b$}} 
 -\mbox{\boldmath{$b$}} \otimes \mbox{\boldmath{$a$}} 
\end{equation*} 
and isomorphic to $\bigwedge^2 E^4$ in Section~\ref{sect:ocs}. 

Suppose that $V$ is oriented and has a fixed inner product $h$. 
Let $(\xi_1 , \xi_2 , \xi_3 , \xi_4 )$ be 
an ordered orthonormal basis of $(V, h)$ giving the orientation. 
Let $\bigwedge^2_{\pm} V$ be the subspaces of $\bigwedge^2 V$ 
generated by 
\begin{equation} 
\begin{split} 
  \Omega_{\pm , 1} & 
:=\dfrac{1}{\sqrt{2}} (\xi_1 \wedge \xi_2 \pm \xi_3 \wedge \xi_4 ), \\ 
  \Omega_{\pm , 2} & 
:=\dfrac{1}{\sqrt{2}} (\xi_1 \wedge \xi_3 \pm \xi_4 \wedge \xi_2 ), \\ 
  \Omega_{\pm , 3} & 
:=\dfrac{1}{\sqrt{2}} (\xi_1 \wedge \xi_4 \pm \xi_2 \wedge \xi_3 ) 
\end{split} 
\label{xiwedgexi} 
\end{equation} 
respectively. 
Then $\bigwedge^2_{\pm} V$ do not depend on 
the choice of $(\xi_1 , \xi_2 , \xi_3 , \xi_4 )$, 
which can be explained by 
the double covering $SO(4)\longrightarrow SO(3)\times SO(3)$. 
We have the decomposition $\bigwedge^2 V 
                          =\bigwedge^2_+ V \oplus \bigwedge^2_- V$. 
The inner product $h$ of $V$ induces 
an  inner product $\hat{h}$ of $\bigwedge^2 V$ 
such that the six elements in \eqref{xiwedgexi} form 
an orthonormal basis of $\bigwedge^2 V$. 

Let $N$ be an oriented $4$-dimensional Riemannian manifold 
and let $h$ be its Riemannian metric. 
Let $T\!N$ be the tangent bundle of $N$. 
Then we can define (the bundle of) 
the two-fold exterior power $\bigwedge^2 T\!N$ of $T\!N$, 
considering the two-fold exterior power of the tangent space of $N$ 
at each point. 
We have the bundle decomposition $\bigwedge^2 T\!N 
                                 =\bigwedge^2_+ T\!N \oplus 
                                  \bigwedge^2_- T\!N$. 
The \textit{twistor spaces\/} associated with $(N, h)$ are 
the unit sphere bundles of the subbundles $\bigwedge^2_{\pm} T\!N$. 
Sections of the twistor spaces correspond to almost complex structures 
on $N$ preserving the metric $h$. 
We can refer to \cite{AHS}, \cite{ES} for the twistor spaces 
associated with oriented $4$-dimensional Riemannian manifolds. 

Let $M$ be a Riemann surface and 
let $F:M\longrightarrow N$ be a conformal immersion of $M$ into $N$. 
Then the \textit{twistor spaces\/} associated with 
the pullback bundle $F^*\!T\!N$ are the unit sphere bundles 
of the subbundles $\bigwedge^2_{\pm} F^*\!T\!N$ 
of $\bigwedge^2 F^*\!T\!N$. 
The \textit{twistor lifts\/} $F_{\pm}$ of $F$ are defined as 
sections of the twistor spaces associated with $F^*\!T\!N$. 
The Levi-Civita connection $\nabla$ of $h$ induces 
a connection $\hat{\nabla}$ of $\bigwedge^2 F^*\!T\!N$. 
In addition, $\hat{\nabla}$ gives 
connections of $\bigwedge^2_{\pm} F^*\!T\!N$. 
Suppose that $F$ is minimal. 
Then an analogous result to Theorem~\ref{thm:F} holds, 
based on \cite{friedrich} (see \cite{ando(6)} for details). 
In this result, 
the condition corresponding to the condition (e) in Theorem~\ref{thm:F} is 
stated in terms of the horizontality of the twistor lifts. 
See \cite{bryant} in the case of $N=S^4$: 
each twistor space associated with $S^4$ is considered to be 
the homogeneous space $SO(5)/U(2)$, 
which is rewritten into $Sp(2)/U(2)=\mbox{\boldmath{$C$}}\!P^3$ 
by the double covering $Sp(2)\longrightarrow SO(5)$; 
the composition of each twistor lift of $F:M\longrightarrow S^4$ 
and the Penrose twistor map 
$T: \mbox{\boldmath{$C$}}\!P^3 \longrightarrow S^4$ is just 
the immersion $F$. 
Suppose that $N$ is hyperK\"{a}hler: 
for example, $K3$-surfaces admit hyperK\"{a}hler structures. 
Then one of the twistor spaces associated with the pullback bundle 
by $F:M\longrightarrow N$ is the product bundle $M\times S^2$, 
and the corresponding twistor lift is a holomorphic map 
from $M$ into $\mbox{\boldmath{$C$}}\!P^1$ (see \cite{ES}, \cite{ando(7)}). 
Therefore Theorem~\ref{thm:HO} is the result in the special case. 
In addition, $F$ is an isotropic minimal immersion 
compatible with the orientation given by the hyperK\"{a}hler structure 
if and only if $F$ is a complex curve 
with respect to an orthogonal complex structure obtained from 
the hyperK\"{a}hler structure (\cite{friedrich}, \cite{ando(2)}). 
There exist a hyperK\"{a}hler $4$-manifold and an isotropic minimal surface 
in it which is not a complex curve with respect to 
any orthogonal complex structure 
obtained from the hyperK\"{a}hler structure (see \cite{AH}, \cite{MW}). 
Suppose that $N$ is a K\"{a}hler surface. 
Then $F$ is a complex curve if and only if 
$F$ is an isotropic minimal immersion 
which is compatible with the orientation given by the K\"{a}hler structure 
and has at least one complex point (\cite{ando(2)}). 

\subsection{Neutral cases} 

Let $V$ be a $4$-dimensional vector space  
and $h$ a nondegenerate symmetric bilinear function on $V$. 
We call $h$ a \textit{neutral inner product\/} 
if there exists a basis $\{ \xi_1 , \xi_2 , \xi_3 , \xi_4 \}$ of $V$ 
satisfying 
\begin{equation} 
 h(\xi_i, \xi_j ) 
=\left\{ \begin{array}{rl} 
          1 & (i=j=1, 2), \\ 
         -1 & (i=j=3, 4), \\ 
          0 & (i\not= j). 
           \end{array} 
 \right. 
\label{nip} 
\end{equation} 
Let $h$ be a neutral inner product of $V$. 
Then such a basis of $V$ as $\{ \xi_1 , \xi_2 , \xi_3 , \xi_4 \}$ 
with \eqref{nip} is called a \textit{pseudo-orthonormal basis\/} 
with respect to $h$. 
Suppose that $V$ is oriented and has a fixed neutral inner product $h$. 
Let $(\xi_1 , \xi_2 , \xi_3 , \xi_4 )$ be 
an ordered pseudo-orthonormal basis of $(V, h)$ giving the orientation. 
Let $\bigwedge^2_{\pm , h} V$ be the subspaces of 
the two-fold exterior power $\bigwedge^2 V$ of $V$ 
generated by $\Omega_{\mp , 1}$, $\Omega_{\pm , 2}$, $\Omega_{\pm , 3}$ 
as in \eqref{xiwedgexi} respectively. 
Then $\bigwedge^2_{\pm , h} V$ do not depend on 
the choice of $(\xi_1 , \xi_2 , \xi_3 , \xi_4 )$, 
which can be explained by 
the double covering $SO_0 (2, 2)\longrightarrow SO_0(1, 2)\times SO_0(1, 2)$. 
We have the decomposition $\bigwedge^2 V 
                          =\bigwedge^2_{+, h} V \oplus 
                           \bigwedge^2_{-, h} V$. 
The neutral inner product $h$ of $V$ induces 
a nondegenerate symmetric bilinear function $\hat{h}$ on $\bigwedge^2 V$ 
satisfying 
\begin{itemize} 
\item{$\bigwedge^2_{\pm , h} V$ are orthogonal to each other 
with respect to $\hat{h}$, that is, 
$\hat{h} (\Omega_+ , \Omega_- )=0$ 
for $\Omega_{\pm} \in \bigwedge^2_{\pm , h} V$;} 
\item{$\Omega_{\mp , 1}$, $\Omega_{\pm , 2}$, $\Omega_{\pm , 3}$ form 
pseudo-orthonormal bases of $\bigwedge^2_{\pm , h} V$ respectively 
with respect to the restrictions of $\hat{h}$, that is, 
for $\varepsilon$, $\varepsilon' \in \{ +, -\}$ and $i\not= j$, 
\begin{equation*} 
\hat{h} (\Omega_{\varepsilon , 1} , \Omega_{\varepsilon  , 1} )= 1, \quad 
\hat{h} (\Omega_{\varepsilon , 2} , \Omega_{\varepsilon  , 2} )= 
\hat{h} (\Omega_{\varepsilon , 3} , \Omega_{\varepsilon  , 3} )=-1, \quad 
\hat{h} (\Omega_{\varepsilon , i} , \Omega_{\varepsilon' , j} )= 0. 
\end{equation*}} 
\end{itemize} 
We say that $\hat{h}$ is 
an \textit{indefinite inner product of $\bigwedge^2 V$ 
with signature\/} $(2, 4)$. 
Then the restrictions of $\hat{h}$ on $\bigwedge^2_{\pm , h} V$ are 
\textit{indefinite inner products with signature\/} $(1, 2)$. 

Let $N$ be a $4$-dimensional manifold 
and $h$ a nondegenerate symmetric $2$-tensor field on $N$. 
We call $h$ a \textit{neutral metric\/} 
if for each point $q$ of $N$, 
$h$ gives a neutral inner product of the tangent space $T_q N$ at $q$. 
Let $h$ be a neutral metric on $N$. 
Then we call $(N, h)$ a \textit{neutral manifold}. 
A simplest example of a $4$-dimensional neutral manifold is 
the neutral Euclidean $4$-space. 
We define a nondegenerate bilinear function on $\mbox{\boldmath{$R$}}^4$ 
by 
\begin{equation} 
  \langle \mbox{\boldmath{$x$}} , \mbox{\boldmath{$y$}} \rangle_{2, 2} 
:=x^1 y^1 +x^2 y^2 -x^3 y^3 -x^4 y^4 
\label{nmetE} 
\end{equation} 
for $\mbox{\boldmath{$x$}} =(x^1 , x^2 , x^3 , x^4 )$, 
    $\mbox{\boldmath{$y$}} =(y^1 , y^2 , y^3 , y^4 ) 
     \in \mbox{\boldmath{$R$}}^4$. 
We call the bilinear function $\langle \ , \ \rangle_{2, 2}$ 
defined in \eqref{nmetE} 
the \textit{standard neutral inner product\/} 
of $\mbox{\boldmath{$R$}}^4$, and 
we call $\mbox{\boldmath{$R$}}^4$ 
equipped with $\langle \ , \ \rangle_{2, 2}$ 
the \textit{neutral Euclidean $4$-space}, 
which is denoted by $E^4_2$. 
The standard neutral inner product $\langle \ , \ \rangle_{2, 2}$ 
of $E^4_2$ defines a unique neutral metric 
on $\mbox{\boldmath{$R$}}^4$ as a $4$-dimensional manifold 
by parallelism, 
which is called the \textit{standard neutral metric\/} on $E^4_2$. 
Hence we consider $E^4_2$ to be a $4$-dimensional neutral manifold. 
Clearly, $E^4_2$ is considered to be 
the neutral analogue of Euclidean $4$-space $E^4$. 

Let $N$ be an oriented $4$-dimensional neutral manifold 
and let $h$ be its neutral metric. 
The two-fold exterior power $\bigwedge^2 T\!N$ of $T\!N$ is 
decomposed into 
\begin{equation} 
 \textstyle\bigwedge^2 T\!N 
=\textstyle\bigwedge^2_{+, h} T\!N\oplus 
 \textstyle\bigwedge^2_{-, h} T\!N. 
\label{2epdecomp} 
\end{equation} 
The \textit{space-like twistor spaces\/} associated with $(N, h)$ are 
fiber bundles in $\bigwedge^2_{\pm , h} T\!N$ given by 
\begin{equation} 
 U_+ (\textstyle\bigwedge^2_{\pm , h} T\!N) 
=\{ \Omega \in \textstyle\bigwedge^2_{\pm , h} T\!N \ | \ 
    \hat{h} (\Omega , \Omega )=1\} 
\label{sts} 
\end{equation} 
respectively. 
A fiber of each space-like twistor space is 
a hyperboloid of two sheets. 
Sections of the space-like twistor spaces correspond to 
almost complex structures on $N$ preserving the neutral metric $h$. 
We can refer to \cite{BDM} for the space-like twistor spaces 
(hyperbolic twistor spaces in \cite{BDM}). 
 
Let $M$ be a connected Riemann surface and 
let $F:M\longrightarrow N$ be a space-like and conformal immersion 
of $M$ into $N$. 
Then the \textit{space-like twistor spaces\/} associated with 
the pullback bundle $F^*\!T\!N$ are 
fiber bundles in $\bigwedge^2_{\pm} F^*\!T\!N$ 
defined in a similar way to \eqref{sts}. 
The \textit{space-like twistor lifts\/} $F_{\pm}$ of $F$ are defined as 
sections of the space-like twistor spaces associated with $F^*\!T\!N$. 
Suppose that $F$ has zero mean curvature vector. 
Then an analogous result to Theorem~\ref{thm:F}  
holds (see \cite{ando(4)}, \cite{ando(6)}). 
Suppose that $N$ is neutral hyperK\"{a}hler: 
see \cite{DGMY}, \cite{kamada} 
for neutral hyperK\"{a}hler $4$-manifolds. 
Then one of the space-like twistor spaces associated with the pullback bundle 
by $F:M\longrightarrow N$ is a product bundle, 
and the corresponding space-like twistor lift is a holomorphic map 
from $M$ into 
the one-dimensional complex hyperbolic space $\mbox{\boldmath{$C$}}\!H^1$ 
(\cite{ando(7)}), 
which yields an analogue of Theorem~\ref{thm:HO} 
in the neutral Euclidean $4$-space $E^4_2$. 
See \cite{ando(4)} for complex curves 
in neutral hyperK\"{a}hler $4$-manifolds and neutral K\"{a}hler surfaces. 

Let $(N, h)$ be an oriented $4$-dimensional neutral manifold. 
Related to the decomposition \eqref{2epdecomp}, 
the \textit{time-like twistor spaces\/} associated with $(N, h)$ are 
fiber bundles in $\bigwedge^2_{\pm , h} T\!N$ given by 
\begin{equation} 
 U_- (\textstyle\bigwedge^2_{\pm , h} T\!N) 
=\{ \Omega \in \textstyle\bigwedge^2_{\pm , h} T\!N \ | \ 
    \hat{h} (\Omega , \Omega )=-1\} 
\label{tts} 
\end{equation} 
respectively. 
A fiber of each time-like twistor space is 
a hyperboloid of one sheet. 
Sections of the time-like twistor spaces correspond to 
almost paracomplex structures on $N$ reversing the neutral metric $h$, 
and see \cite{CFG} for paracomplex structures. 
We can refer to \cite{HM}, \cite{JR} for the time-like twistor spaces. 

Let $M$ be a Lorentz surface. 
Lorentz surfaces are analogues of Riemann surfaces, and 
when we consider Lorentz surfaces instead of Riemann surfaces, 
we need to replace 
the complex plane $\mbox{\boldmath{$C$}}$ and holomorphicity 
with the paracomplex plane $\check{\mbox{\boldmath{$C$}}}$ 
and paraholomorphicity respectively. 
See the appendix of \cite{ando(4)} for Lorentz surfaces, 
and see \cite{ando(7)} 
for a hyperboloid of one sheet $\check{\mbox{\boldmath{$C$}}}\!H^1$ 
as a Lorentz surface, 
which is an analogue of the Riemann sphere $\mbox{\boldmath{$C$}}\!P^1$. 
Let $F:M\longrightarrow N$ be a time-like and conformal immersion 
of $M$ into $N$. 
Then the \textit{time-like twistor spaces\/} associated with 
the pullback bundle $F^*\!T\!N$ are 
fiber bundles in $\bigwedge^2_{\pm} F^*\!T\!N$ 
defined in a similar way to \eqref{tts}. 
The \textit{time-like twistor lifts\/} $F_{\pm}$ of $F$ are defined as 
sections of the time-like twistor spaces associated with $F^*\!T\!N$. 
Suppose that $F$ has zero mean curvature vector. 
Then an analogous result to Theorem~\ref{thm:F} does not hold: 
see \cite{ando(4)}, \cite{ando(6)}, \cite{ando(7)}, and 
see \cite{ando(9)} for characterizations of 
a time-like and conformal immersion with zero mean curvature vector 
into a $4$-dimensional neutral space form 
such that the covariant derivatives of the two time-like twistor lifts 
are zero or light-like. 
Suppose that $N$ is neutral hyperK\"{a}hler. 
Then one of the time-like twistor spaces associated with the pullback bundle 
by $F:M\longrightarrow N$ is a product bundle, 
and the corresponding time-like twistor lift is a paraholomorphic map 
from $M$ into $\check{\mbox{\boldmath{$C$}}}\!H^1$ (\cite{ando(7)}), 
which yields an analogue of Theorem~\ref{thm:HO} in $E^4_2$. 
See \cite{ando(4)} for paracomplex curves 
in neutral hyperK\"{a}hler $4$-manifolds and paraK\"{a}hler surfaces. 

\subsection{Lorentzian cases} 

Let $V$, $h$ be as in the beginning of the previous subsection. 
We call $h$ a \textit{Lorentz-Minkowski inner product\/} 
if there exists a basis $\{ \xi_1 , \xi_2 , \xi_3 , \xi_4 \}$ of $V$ 
satisfying 
\begin{equation} 
 h(\xi_i, \xi_j ) 
=\left\{ \begin{array}{rl} 
          1 & (i=j=1, 2, 3), \\ 
         -1 & (i=j=4), \\ 
          0 & (i\not= j). 
           \end{array} 
 \right. 
\label{LMip} 
\end{equation} 
Let $h$ be a Lorentz-Minkowski inner product of $V$. 
Then such a basis of $V$ as $\{ \xi_1 , \xi_2 , \xi_3 , \xi_4 \}$ 
with \eqref{LMip} is called a \textit{pseudo-orthonormal basis\/} 
with respect to $h$, as in the previous subsection. 
Suppose that $V$ is oriented 
and has a fixed Lorentz-Minkowski inner product $h$. 
Let $(\xi_1 , \xi_2 , \xi_3 , \xi_4 )$ be 
an ordered pseudo-orthonormal basis of $(V, h)$ giving the orientation. 
The Lorentz-Minkowski inner product $h$ of $V$ induces 
a nondegenerate symmetric bilinear function $\hat{h}$ on $\bigwedge^2 V$ 
such that $\xi_i \wedge \xi_j$ ($1\leq i<j\leq 4$) are orthogonal to 
one another so that 
\begin{equation*} 
 \hat{h} (\xi_i \wedge \xi_j , \xi_i \wedge \xi_j ) 
=\left\{ \begin{array}{rl} 
          1 & (1\leq i<j\leq 3), \\ 
         -1 & (j=4). 
           \end{array} 
 \right. 
\end{equation*} 
We say that $\hat{h}$ is 
a  \textit{neutral inner product\/} or 
an \textit{indefinite inner product with signature\/} 
$(3, 3)$ of $\bigwedge^2 V$. 
We set 
\begin{equation*} 
   L(\textstyle\bigwedge^2 V) 
:=\{ \Omega \in \textstyle\bigwedge^2 V\setminus \{ 0\} \ | \ 
     \hat{h} (\Omega , \Omega )=0\} . 
\end{equation*} 
Then $L(\bigwedge^2 V)$ is the set of light-like elements 
of $\bigwedge^2 V$ with respect to $\hat{h}$, 
which is a hypersurface in $\bigwedge^2 V$.  
Let $\mathcal{B} (V, h)$ be 
the set of ordered pseudo-orthonormal bases of $(V, h)$ 
giving the orientation. 
Then for $(\xi_1 , \xi_2 , \xi_3 , \xi_4 )\in \mathcal{B} (V, h)$, 
$\Omega_{\pm , 1}$ as in \eqref{xiwedgexi} is 
an element of $L(\bigwedge^2 V)$. 
We set 
\begin{equation} 
\mathcal{L}_{\pm} (\textstyle\bigwedge^2 V) 
:=\{ \Omega_{\pm , 1} \ | \ 
    (\xi_1 , \xi_2 , \xi_3 , \xi_4 )\in \mathcal{B} (V, h)\} . 
\label{Lpm} 
\end{equation} 
Then $\mathcal{L}_{\pm} (\bigwedge^2 V)$ are 
hypersurfaces in $L(\bigwedge^2 V)$ 
which are $SO(3, 1)$-orbits (see \cite{ando(10)} 
for $SO(3, 1)$-orbits in $L(\bigwedge^2 V)$), 
and $4$-dimensional neutral manifolds. 

Let $N$ be a $4$-dimensional manifold 
and $h$ a nondegenerate symmetric $2$-tensor field on $N$. 
We call $h$ a \textit{Lorentzian metric\/} 
if for each point $q$ of $N$, 
$h$ gives a Lorentz-Minkowski inner product of the tangent space $T_q N$ 
at $q$. 
Let $h$ be a Lorentzian metric on $N$. 
Then we call $(N, h)$ a \textit{Lorentzian manifold}. 
A simplest example of a $4$-dimensional Lorentzian manifold is 
the Lorentz-Minkowski $4$-space. 
We define a nondegenerate bilinear function on $\mbox{\boldmath{$R$}}^4$ 
by 
\begin{equation} 
  \langle \mbox{\boldmath{$x$}} , \mbox{\boldmath{$y$}} \rangle_{3, 1} 
:=x^1 y^1 +x^2 y^2 +x^3 y^3 -x^4 y^4 
\label{LMmetE} 
\end{equation} 
for $\mbox{\boldmath{$x$}} =(x^1 , x^2 , x^3 , x^4 )$, 
    $\mbox{\boldmath{$y$}} =(y^1 , y^2 , y^3 , y^4 ) 
     \in \mbox{\boldmath{$R$}}^4$. 
We call the bilinear function $\langle \ , \ \rangle_{3, 1}$ 
defined in \eqref{LMmetE} 
the (\textit{standard\/}) \textit{Lorentz-Minkowski inner product\/} 
of $\mbox{\boldmath{$R$}}^4$, and 
we call $\mbox{\boldmath{$R$}}^4$ 
equipped with $\langle \ , \ \rangle_{3, 1}$ 
the \textit{Lorentz-Minkowski $4$-space}, 
which is denoted by $E^4_1$. 
The Lorentz-Minkowski inner product $\langle \ , \ \rangle_{3, 1}$ 
of $E^4_1$ defines a unique Lorentzian metric 
on $\mbox{\boldmath{$R$}}^4$ as a $4$-dimensional manifold 
by parallelism, 
which is called the (\textit{standard\/}) \textit{Lorentz-Minkowski metric\/} 
on $E^4_1$. 
Hence we consider $E^4_1$ to be a $4$-dimensional Lorentzian manifold, 
and $E^4_1$ is considered to be 
the Lorentzian analogue of Euclidean $4$-space $E^4$. 

Let $(N, h)$ be a $4$-dimensional Lorentzian manifold. 
Let $M$ be a Riemann surface and 
let $F:M\longrightarrow N$ be a space-like and conformal immersion 
of $M$ into $N$. 
Let $\mathcal{L}_{\pm} (\bigwedge^2 F^*\!T\!N)$ be 
fiber bundles in $\bigwedge^2 F^*\!T\!N$ 
such that the fibers on each point $p$ of $M$ are given as in \eqref{Lpm} 
with $V=F^*\!T_p N$. 
The \textit{lifts\/} $F_{\pm}$ of $F$ are defined as 
sections of $\mathcal{L}_{\pm} (\bigwedge^2 F^*\!T\!N)$ respectively. 
Suppose that $F$ has zero mean curvature vector. 
If $N=E^4_1$, 
then the lifts $F_{\pm}$ are considered 
to be maps from $M$ into $\mathcal{L}_{\pm} (\bigwedge^2 E^4_1 )$, 
and in addition, 
they are holomorphic with respect to parallel almost complex structures 
on $\mathcal{L}_{\pm} (\bigwedge^2 E^4_1 )$ (\cite{ando(5)}, \cite{ando(8)}), 
which can be considered to be 
a result in $E^4_1$ analogous to Theorem~\ref{thm:HO}. 
See \cite{ando(3)}, \cite{ando(5)}, \cite{AnHa} for characterizations of 
a space-like and conformal immersion with zero mean curvature vector 
into a $4$-dimensional Lorentzian space form 
such that the covariant derivatives of the two lifts are zero or light-like, 
and perpendicular to each other. 

Let $M$ be a Lorentz surface and 
let $F:M\longrightarrow N$ be a time-like and conformal immersion 
of $M$ into $N$. 
Then, as in the space-like case, 
the \textit{lifts\/} $F_{\pm}$ of $F$ are defined as 
sections of $\mathcal{L}_{\pm} (\bigwedge^2 F^*\!T\!N)$ respectively. 
Suppose that $F$ has zero mean curvature vector. 
If $N=E^4_1$, 
then the lifts $F_{\pm}$ are paraholomorphic 
with respect to parallel almost paracomplex structures 
on $\mathcal{L}_{\pm} (\bigwedge^2 E^4_1 )$ (\cite{ando(5)}, \cite{ando(8)}), 
which is the time-like analogue of the corresponding result 
in the previous paragraph. 
Refer to \cite{ando(5)}, \cite{AnHa} for characterizations of 
a time-like and conformal immersion with zero mean curvature vector 
into a $4$-dimensional Lorentzian space form 
such that the covariant derivatives of the two lifts are zero or light-like, 
and perpendicular to each other.

\vspace{4mm} 

\par\noindent 
\footnotesize{Faculty of Advanced Science and Technology, 
              Kumamoto University \\ 
              2--39--1 Kurokami, Chuo-ku, Kumamoto 860--8555 Japan} 

\par\noindent  
\footnotesize{E-mail address: andonaoya@kumamoto-u.ac.jp} 

\end{document}